\documentclass{amsart}
\usepackage{amsmath,amssymb}
\usepackage{tikz}

\usepackage{hyperref}
\hypersetup{colorlinks=true,linkcolor=red, citecolor=green, urlcolor=cyan}

\newtheorem{theorem}{Theorem}[section]
\newtheorem{prop}[theorem]{Proposition}

\newtheorem{lemma}[theorem]{Lemma}
\theoremstyle{definition}
\newtheorem{defi}[theorem]{Definition}
\theoremstyle{remark}

\numberwithin{equation}{section}

\DeclareMathOperator{\Id}{Id}
\DeclareMathOperator{\supp}{supp}
\DeclareMathOperator{\tr}{tr}
\DeclareMathOperator{\diag}{diag}
\DeclareMathOperator{\spn}{span}
\DeclareMathOperator*{\wlim}{w-lim}
\DeclareMathOperator{\dist}{dist}

\DeclareMathOperator{\CS}{\Pi}

\newcommand{\be}{\begin{equation}}
\newcommand{\ee}{\end{equation}}
\newcommand{\bes}{\begin{equation*}}
\newcommand{\ees}{\end{equation*}}

\newcommand{\R}{\mathbb{R}} 
\newcommand{\C}{\mathbb{C}} 
\newcommand{\N}{\mathbb{N}}

\newcommand{\Z}{\mathbb{Z}}

\newcommand{\p}{\partial}
\newcommand{\eps}{\varepsilon}
\newcommand{\la}{\lambda}
\newcommand{\pv}{\textrm{p.v.}}
\newcommand{\vk}{\varkappa}

\newcommand{\va}[1]{\bigl\lvert#1\bigr\rvert} 
\newcommand{\ps}[2]{  \bigl\langle #1\,,\, #2  \bigr\rangle} 
\newcommand{\qtq}[1]{\quad\text{#1}\quad}

\newcommand{\Hohp}{H^{\smash[t]{1/2}}_+} 

\let\Re=\undefined\DeclareMathOperator{\Re}{Re}
\let\Im=\undefined\DeclareMathOperator{\Im}{Im}

\newcounter{smalllist}
\newenvironment{SL}{\begin{list}{\hss\upshape(\roman{smalllist})\hss}{%
\setlength{\topsep}{0mm}\setlength{\parsep}{0mm}\setlength{\itemsep}{0mm}%
\setlength{\labelwidth}{1.75em}\setlength{\labelsep}{\the\fontdimen2\font}\setlength{\leftmargin}{\the\labelwidth}\addtolength{\leftmargin}{\labelsep}%
\setlength{\itemindent}{0em}\usecounter{smalllist}%
}}{\end{list}}

\allowdisplaybreaks

\begin{document}

\title{Orbital stability of Benjamin--Ono multisolitons}

\author{Rana Badreddine}
\address{Department of Mathematics, University of California, Los Angeles, CA 90095, USA}
\email{badreddine@math.ucla.edu}

\author{Rowan Killip}
\address{CEREMADE, CNRS, Universit\'e Paris Dauphine--PSL, Place du Mar\'echal de Lattre de Tassigny, 75016 Paris, France \&  Department of Mathematics, University of California, Los Angeles, CA 90095, USA}
\email{killip@ceremade.dauphine.fr}

\author{Monica Vi\c{s}an}
\address{Department of Mathematics, University of California, Los Angeles, CA 90095, USA}
\email{visan@math.ucla.edu}

\begin{abstract}
We show that multisoliton solutions to the Benjamin--Ono equation are uniformly orbitally stable in $H^s(\R)$ for every $-\tfrac12<s\leq \frac12$.  This improves the regularity required for stability up to the sharp well-posedness threshold; previous work (even on single solitons) had required $s\geq \frac12$.

One key ingredient in our argument is a new variational characterization of multisolitons.  A second ingredient is the extension to low-regularity slowly-decaying solutions of the Wu identity on eigenfunctions of the Lax operator. This extension also allows us to clarify the spectral type of the Lax operator for such potentials, by precluding embedded eigenvalues.
\end{abstract}
\maketitle


\section{Introduction}\label{S:1}
We consider the Benjamin--Ono equation
\begin{equation}\label{BO}\tag{BO}
    \partial_t u = H\partial_x^2 u - 2u\partial_x u, \qquad u: \R_t\times\R_x\to\R, 
\end{equation}
which was introduced in \cite{Benjamin1967,Davis1967} to  describe the motion of internal waves in stratified fluids of great total depth. The function $u$ represents the elevation of the interface relative to equilibrium. The operator $H$ is the Hilbert transform, defined via
\begin{equation}\label{HT1}
    Hf(x) = \pv \frac{1}{\pi} \int \frac{f(y)}{x-y}\,dy 
\end{equation}
and expressed in Fourier variables \eqref{FT} via
\begin{equation}\label{HT2}
    \widehat{Hf}(\xi) = -i\, \text{sgn}(\xi)\widehat{f}(\xi).
\end{equation}

Both Benjamin \cite{Benjamin1967} and Davis--Acrivos \cite{Davis1967} carried out laboratory experiments to corroborate the predictions of \eqref{BO} and witnessed the solitary wave solutions that they had predicted mathematically.  These have the form
\begin{equation}\label{soliton Qc}
Q_{\lambda,c}(t,x) = \frac{-4\lambda}{1+ 4\lambda^2(x-c+2\lambda t)^2} \qtq{with} \lambda<0 \qtq{and} c\in\R.
\end{equation}
We see that $Q_{\lambda,c}$ has a Lorentzian profile of width $|\lambda|^{-1}$. It moves to the right with speed $2|\lambda|$ and $c$ marks the center of the wavepacket at time $t=0$.  Our decision to use $\lambda<0$ as a parameter will be justified later --- it is the eigenvalue of the corresponding Lax operator!

These researchers drew particular attention to the solitary wave's ``remarkable property of persistence'' \cite[p.~564]{Benjamin1967} and ``... the ease with which it could be generated.'' \cite[p.~604]{Davis1967}.  As we will discuss in due course, the stability they observed has received rigorous mathematical justification.  In this paper, we will develop this theory further in two directions: (1) by treating more complex coherent structures and (2) by allowing a wider class of initial perturbations.  Concretely, we will consider general multisoliton solutions to \eqref{BO} and allow much `noisier' perturbations by placing much weaker demands on the high-frequency components of solutions.

Stability is a dynamical property; it can only be discussed for perturbations that lead to an unambiguous evolution.  For this reason, we must work in spaces where the equation is well-posed, that is, solutions to the initial-value problem exist, are unique, and depend continuously on the initial data.  The determination of when canonical dispersive equations are well- or ill-posed has been a major mathematical quest in recent decades.  Works investigating the \eqref{BO} model include  \cite{Abdelouhab1989, Burq2006,  Ginibre1989a, Ginibre1989, Ginibre1991, Ifrim2019, Ionescu2007, Iorio1986, Kenig2003, Koch2003, Koch2005, Molinet2007, Molinet2008, Molinet2012, Molinet2009, Molinet2001, Ponce1991, Saut1979, Talbut2021, Tao2004}.  Nevertheless, it was only very recently that the optimal well-posedness result in the scale of Sobolev spaces $H^s(\R)$ was achieved.  This was obtained in \cite{Killip2024} and is reproduced in Theorem~\ref{T:KLV main} below; the sharp result for periodic initial data was shown earlier in \cite{GKTActa}.

\begin{theorem}[\cite{Killip2024}]\label{T:KLV main}
Fix $s>-\frac{1}{2}$.  The equation~\eqref{BO} is globally well-posed for initial data in $H^s(\R)$.
\end{theorem}

Ill-posedness for $s<-\frac12$ can be deduced from the explicit solution \eqref{soliton Qc}: At time $t=0$, the solution converges in $H^s(\R)$ to a delta function as $\lambda\to-\infty$, but it does not converge at any other time.  Ill-posedness at the threshold $s=-\frac12$ has been shown in the periodic setting \cite{GKTActa}; we believe that the same physical mechanism yields ill-posedness in $H^{-1/2}(\R)$, but this remains to be shown rigorously.

The identification of the optimal well-posedness space sets us a clear target: prove stability in $H^s(\R)$ for any $s>-\frac12$. This is what we will achieve, not only for single solitary waves, but for the broader class of multisolitons that we will now describe.

The observation of solitary water waves dates back to the early 19th century.  By contrast, the name \emph{soliton} was coined in 1965 by Kruskal and Zabusky \cite{ZabuskyKruskal}.  In numerical simulations of the Korteweg--de Vries equation, they observed solitary waves ``preserving their identity through numerous interactions''.   The waves were behaving like particles and the name \emph{soliton} was chosen to reflect this behavior.  This observation was the launching-pad for the subsequent discovery of the complete integrability of the KdV equation.  In particular, KdV was discovered to admit exact solutions that are nonlinear superpositions of finitely many solitons; these are the multisolitons.

In the paper \cite{Ono1975} that stimulated much subsequent work, Ono explicitly asks if the traveling wave solutions \eqref{soliton Qc} to \eqref{BO} are solitons.  He also highlights the slow spatial decay of the wave profile, which distinguishes it sharply from the exponential decay of KdV and NLS solitons. 

The answer to Ono's question came quickly in \cite{Case1979,MR516327,Joseph1977,Matsuno1979,MeissPereira}.  Solitary waves for Benjamin--Ono do indeed interact like solitons and there are corresponding exact multisoliton solutions.  The following description of these multisolitons is adapted from \cite[\S3.1]{MR759718}:

\begin{defi}[Multisoliton solutions]\label{D:multisoliton}
Fix $N\geq 1$.  Given $\la_1< \cdots <\la_N<0$ and $N$ real parameters $c_1, \ldots, c_N\in \R$, let $\Lambda=\{\lambda_1,\ldots,\lambda_N\}$ and
\begin{align}\label{E:Qbc}
Q_{\Lambda, \vec c} \,(x)= - 2\Im \tfrac{d}{dx} \ln \det\bigl[M(x)\bigr] ,
\end{align}
where $M(x)$ is the $N\times N$ matrix with entries
\begin{align}\label{matrix A}
M_{jk}(x) = \begin{cases}
-i(x-c_j) - \tfrac{1}{2\la_j}, \qquad &j=k\\
-\frac{1}{\la_j-\la_k},  \qquad &j\neq k.
\end{cases}
\end{align}
The unique solution to \eqref{BO} with initial data $u(0,x) = Q_{\Lambda, \vec c}\,(x)$ is the $N$-soliton solution
\begin{align}\label{E:Qbc(t)}
u(t,x) = Q_{\Lambda, \vec c(t)}(x) \quad\text{where} \quad c_j(t) = c_j - 2\la_j t.
\end{align}
\end{defi}

The solutions \eqref{E:Qbc(t)} are most easily identified as multisolitons by observing their long-time behaviour.  As $t\to\pm\infty$, the matrix $M$ becomes increasingly dominated by its diagonal and one can verify that the solution is increasingly well approximated, in $L^2(\R)$, by a linear combination of $N$ well-separated solitary waves as defined in \eqref{soliton Qc}.  In fact, these waves have precisely the parameters $\Lambda$ and centers $c_1, \ldots, c_n$.  The choice of $L^2(\R)$ is significant: this norm is conserved by the flow, so convergence in $L^2(\R)$ guarantees that the solution has no radiation component.

For a given set of parameters $\Lambda=\{\lambda_1<\cdots<\lambda_N\}$, varying $\vec c$ leads to a multisoliton manifold
\begin{equation}\label{MSM}
\mathcal M(\Lambda) = \{ Q_{\Lambda, \vec c} : \vec c \in \R^N\},
\end{equation}
which is not only invariant under \eqref{BO}, but under the whole \eqref{BO} hierarchy; see \cite{Sun2021}. 
Our main result, Theorem~\ref{Th: Orbital Stability}, is uniform orbital stability of this manifold: For any $\eps>0$ there is a choice of $\delta>0$ so that solutions with initial data inside a $\delta$-neighborhood of the manifold remain within a distance $\eps$ of the manifold for all time.  Distance will be measured in the $H^s(\R)$ norm.

The word \emph{uniform} captures the fact that the initial data is chosen from a neighborhood of $\mathcal M$ of uniform width $\delta$. This distinction will be useful when we discuss prior work.  The qualifier \emph{orbital} expresses the fact that trajectories remain close to the manifold, rather than to a single solution.  This freedom is necessary, even in the case of a single soliton.  Indeed, a small perturbation of the parameter $\lambda$ in \eqref{soliton Qc} leads to a small change in the initial data; however, even a small modification in $\lambda$ induces a change in speed, and consequently the two trajectories will drift far apart as $t\to\pm\infty$.

\begin{theorem}[Uniform orbital stability of multisolitons]\label{Th: Orbital Stability}
Fix $-\frac12<s\leq \frac12$, an integer $N\geq 1$, and a set $\Lambda$ comprised of negative parameters $\la_1<\cdots<\la_N$. For every \(\eps>0\), there exists $\delta>0$ so that for every initial data $u_0\in H^s(\R)$ satisfying
\[
      \inf_{\vec{c}\in\R^N}\,\|u_0-Q_{\Lambda,\vec{c}}\,\|_{H^s}<\delta,
\]
the corresponding solution \( u(t) \) of \eqref{BO} satisfies
\[
     \sup_{t\in\R}\inf_{\vec{c}\in\R^N}\, \|u(t)-Q_{\Lambda,\vec{c}}\,\|_{H^s}<\eps.
\]
\end{theorem}

The restriction $-\frac12<s\leq \frac12$ is not a fundamental limitation of our method.  Indeed, in Section~\ref{S:HR}, we exhibit two methods for deducing stability at higher regularity from stability at regularities $-\frac12<s<0$. One method is based on equicontinuity, the other on exact conservation laws.   Our decision to treat just this range of $s$ represents a compromise  between generality and readability.  It was important to us to cover the low regularity regime that has been inaccessible to prior methods and, particularly, to reach the sharp well-posedness threshold.  As we will discuss shortly, there are prior results on stability once $s\geq\frac12$.  It was also important for us to treat the particular cases $s=0$ and $s=\frac12$ that align with physically important conservation laws (momentum and energy, respectively).

The proof of Theorem~\ref{Th: Orbital Stability} will be based on the new variational characterization of multisolitons given in Theorem~\ref{Th: variational problem}.  Before presenting an outline of our argument in subsection~\ref{SS:Struct}, we turn to a brief discussion of prior work on orbital stability, with a particular emphasis on the development of the variational approach.

\subsection{Prior work}\label{SS:Prior} Just as solitons were first discovered in experiments modeled by the Korteweg--de Vries equation, so their stability was first investigated in that setting.  Boussinesq \cite{Boussinesq} was already aware of the first three conservation laws for KdV, namely, the conservation of matter, momentum, and energy (to use their modern names).  He referred to the Hamiltonian as the moment of stability because he intuited that solitary waves were minimizers of this, subject to constrained momentum.  In particular, he verified that solitons solve the corresponding Euler--Lagrange equation.  This idea is the root of all subsequent variational approaches to stability, including our own.  

Boussinesq's idea was first made rigorous by Benjamin, working in the KdV setting; see \cite{MR338584,MR386438}.  He proved that KdV solitons are uniformly orbitally stable in the energy space $H^1(\R)$.   He achieved this by showing that solitons are the unique minimizers of the Hamiltonian subject to constrained momentum and by making a quantitative connection between the energy increment above the minimum and the distance to this set of minimizers.  Evidently, this relies on an understanding of the Hessian of the energy at solitons.  The full Hessian is not even positive semi-definite; thus, it is essential that one must only consider variations that are parallel to the constraint manifold.

In their work \cite{MR677997} on nonlinear Schr\"odinger equations, Cazenave--Lions modified this method in a manner that requires only a qualitative form of this trapping property, namely, that all optimizing sequences converge to the soliton manifold.  They employed concentration compactness techniques to overcome the fact the optimizing sequences need not converge.  Indeed, the soliton manifold is not compact, so even sequences of \emph{optimizers} need not converge!


There is a considerable leap in difficulty from treating single solitons to multisolitons.  The pioneers here were Maddocks--Sachs \cite{MR1220540}, who proved that the manifold of $N$-soliton solutions to KdV is uniformly orbitally stable in $H^N(\R)$.  They showed that such multisolitons are local minimizers of the polynomial conservation law for KdV at regularity $H^N(\R)$, constrained by the values of the lower-regularity conservation laws.  (It was shown recently in \cite{MR4892297,MR4612766} that multisolitons are in fact \emph{global} minimizers of this variational problem.)

The paper \cite{KV2022} gave a new variational characterization of KdV multisolitons, based on a renormalization of the reciprocal transmission coefficient.  This characterization, which remains valid in $H^{-1}(\R)$ irrespective of the complexity of the multisoliton, was then used to prove that multisolitons are uniformly orbitally stable in $H^{s}(\R)$ for all $s\geq -1$.  Notice that $s=-1$ corresponds to the optimal space for well-posedness for KdV; see \cite{Killip2019}. The paper \cite{KV2022} is a major inspiration for what is done here.  It also contains a discussion of other approaches to the orbital stability of KdV multisolitons.

Uniform orbital stability of single solitons for \eqref{BO} in the energy space $H^{1/2}(\R)$ was first proved in \cite{MR715035}.  This paper employed the characterization of solitons as minimizers of the energy at constrained momentum and a careful analysis of the corresponding Hessian. This result was revisited in a number of subsequent papers \cite{Albert1992,Albert1999,MR887857,MR2082818,MR897729,MR886343}.  Linearized stability had been shown earlier in \cite{ChenKaup1980}.


The Maddocks--Sachs strategy was applied to \eqref{BO} in \cite{LanWang2025,Matsuno2006,MR1442235,MR2202311}, proving that the manifold associated to any $N$-soliton is uniformly orbitally stable in $H^{N/2}(\R)$.   The key difficulty is understanding the coercivity of the Hessian at every multisoliton.   Once again, it is crucial that only variations preserving the constraints need to be considered.  

One advantage of the variational approach we present in this paper is that we can work at fixed regularity, irrespective of the complexity of the multisoliton.   This virtue is shared by a dynamical approach to orbital stability introduced in \cite{MartelMerle2001,MartelMerleTsai2002} and applied to \eqref{BO} in \cite{GTT2009,KenMart2009}.  This style of argument focuses directly on the large-time regime where the individual solitons that comprise the multisoliton have become well-separated, and will only continue to separate as time progresses.  One strength of this argument is that it has been shown to yield forms of asymptotic stability.  Concretely, \cite{KenMart2009} shows $L^2$ convergence on windows traveling with the solitons for solutions whose initial data is a small $H^{1/2}(\R)$ perturbation of a multisoliton.  One weakness of this argument is that it cannot yield \emph{uniform} orbital stability of $N$-multisolitons because the interval of time needed for the constituent solitons to become favorably configured may be \emph{arbitrarily large} and  is treated via well-posedness theory. 


\subsection{Structure of our argument}\label{SS:Struct}
While the existence and stability of solitary waves is common among dispersive PDE, the observation of solitons (distinguished by their elastic interactions) is a signal of complete integrability.  In our analysis, the complete integrability of \eqref{BO} will enter primarily through the Lax pair representation introduced in 
\cite{Bock1979,Nakamura1979}; our presentation of this material is strongly influenced by the later works \cite{MR3484397} and \cite{Killip2024}.

To describe the Lax operator, we first need some notation.  We write $\CS$ for the Cauchy--Szeg\H{o} projection, that is, the orthogonal projection of $L^2(\R)$ onto the Hardy space
$$
L^2_+(\mathbb{R}) := \{ f \in L^2(\mathbb{R}) :\, \supp(\widehat{\mkern0.3mu f}\mkern1mu) \subseteq [0,+\infty) \}.
$$
Equivalently, $\CS=\tfrac12[{1+iH}]$, where $H$ is the Hilbert transform; see \eqref{HT1} and \eqref{HT2}.  Throughout this paper, we will often use the shorthand
\begin{equation}
u_+ = \CS u .
\end{equation}

For smooth functions, say $u\in H^\infty(\R)$, it is not difficult to verify that
\begin{equation}\label{Lax op}
    L_u f := -i\partial_x f - \Pi(uf)\, \quad \text{and} \quad P_u f := -i\partial_x^2 f - 2\partial_x \Pi(uf) + 2 \Pi(u') f
\end{equation}
define (unbounded) self-adjoint and skew-adjoint operators, respectively, on the Hilbert space $L^2_+(\mathbb{R})$.  Moreover,
\begin{equation}
\text{$u(t)$ is a solution to \eqref{BO}} \iff \tfrac{d}{dt} L_{u(t)} = [P_{u(t)}, L_{u(t)}].  \label{Lax eq}
\end{equation}
We refer to $P_u$ and $L_u$ as the Peter and Lax operators, respectively. 

As shown in \cite{Killip2024}, for general $u\in H^{s}(\R)$, $L_u$ can be realized as a self-adjoint operator that is bounded from below, provided only that $s>-\frac12$; see Theorem~\ref{T:Lu} for a recap of results from \cite{Killip2024} concerning $L_u$.  In particular, these results guarantee that
\begin{equation}\label{beta defn}
\beta\bigl(\kappa; u(t)\bigr) :=\bigl\langle \CS u(t), \bigl(L_{u(t)}+\kappa\bigr)^{-1} \CS u(t) \bigr\rangle
\end{equation}
is well defined for $u\in H^s(\R)$, $s>-\frac12$, and $\kappa$ sufficiently large.  Concretely, when $-\frac12<s\leq0$, it suffices to consider
\begin{align}\label{k0 intro}
  \kappa \geq \kappa_0(u):= C_s \bigl( 1 + \|u\|_{H^s} \bigr)^{\frac2{1+2s}}   .
\end{align}
For $s\geq 0$, the corresponding choice is $\kappa_0(u):= C_s ( 1 + \|u\|_{H^s})^{2} $.  This difference is the first of several inconveniences that would appear if we were to treat all values of $s$ simultaneously.  We have found it most convenient to first restrict attention to $-\frac12<s<0$. We will then treat the remaining values of $s$ in Section~\ref{S:HR}.

A remarkable fact about \eqref{BO} is that $\beta(\kappa; u(t))$ is conserved under the flow.  This was first observed in \cite{MR690743} (in the context of smooth solutions).  As emphasized in \cite{KLV2025}, it can be traced to a very peculiar property of the Peter operator, namely,
\begin{equation}\label{OTP}
\tfrac{d}{dt} \CS u(t) = P_{u(t)} \CS  u(t).
\end{equation}

As discussed in \cite{Kaup1998}, the functional $\beta(\kappa;u)$ serves as a generating function for the traditional sequence of conservation laws associated to \eqref{BO}.  Concretely,
\begin{align}\label{beta expand}
\beta(\kappa; u)= \sum_{n\geq 0} (-1)^n\kappa^{-n-1} \mathcal E_n(u) \qtq{where} \mathcal E_n(u)=\langle u_+,L_u^n\, u_+\rangle.
\end{align}
For our analysis, it is important to control low regularity norms.  This will be discussed in Propositions~\ref{P:beta} and~\ref{P:s pos}.

The generating function $\beta(\kappa;u)$ sits at the heart of our variational characterization of multisolitons.  As noted earlier, solitons are intimately connected with discrete eigenvalues of $L_u$; however, it is not immediately clear that discrete eigenvalues actually contribute to $\beta(\kappa;u)$.  One must fear that $u_+$ is orthogonal to the corresponding eigenfunctions!   This fear is ultimately illusory:

\begin{theorem}[The Wu relations]\label{th: Wu relation}
Fix $s>-\tfrac{1}{2}$, $u\in H^s(\R)$ and let $f$ be an eigenfunction of $L_u$ with eigenvalue $\lambda$.  Then $\lambda<0$,
\be \label{Wu relation Hs}
|\langle u, f \rangle|^2 = -2\pi \lambda \|f\|_{L^2}^2,
\ee
and 
\begin{align}\label{Wu 2}
\sqrt{{2\pi}} \lambda \widehat f(0)  +\langle u, f\rangle =0.
\end{align}
In particular, all eigenvalues of $L_u$ are necessarily negative and simple.
\end{theorem}

In referring to \eqref{Wu relation Hs} and \eqref{Wu 2} as the Wu relations, we are honoring \cite{MR3484397} where this is proved under the assumptions $u\in L^\infty(\R)$ and $\langle x\rangle u\in L^2(\R)$.  Later, Sun \cite{Sun2021} removed the condition $u\in L^\infty(\R)$.

Theorem~\ref{th: Wu relation} constitutes progress on the spectral-theoretic problems posed in \cite[\S1.2]{Killip2024}. In particular, it precludes eigenvalues embedded in the essential spectrum of $L_u$ if $u\in H^s(\R)$ with $s>-\frac12$. By contrast, positive eigenvalues can occur for Schr\"odinger operators already for potentials with $\langle x \rangle^{-1}$ decay.

The proof of Theorem~\ref{th: Wu relation} occupies Section~\ref{S:Wu}; it is essentially self contained, relying only on basic properties of the Lax operator that are recounted in Section~\ref{S:Prelims}. 

Let us turn now to our variational characterization of multisolitons:

\begin{theorem}[Variational characterization of multisolitons]\label{Th: variational problem}
Fix $-\frac12<s<0$, $N\geq 1$, and distinct parameters  \( \la_1<\cdots <\la_N<0\). If $u\in H^s(\R)$ is such that $\lambda_1, \ldots, \lambda_N$ are eigenvalues of the Lax operator $L_u$,
then 
\begin{align}\label{1:45}
\beta(\kappa;u)\geq \sum_{n=1}^N\frac{2\pi\va{\la_n}}{\la_n+\kappa}
\end{align}
whenever $\kappa$ satisfies \eqref{k0 intro}. Moreover, equality holds in \eqref{1:45} for some such $\kappa$ if and only if $u\in \{Q_{\Lambda, \vec c} :\, \vec c\in \R^N\}$ with $\Lambda=\{\lambda_1,\ldots,\lambda_N\}$.
\end{theorem}

This result is proved in Section~\ref{S:VChar}.  The lower bound \eqref{1:45} follows rather directly from the Wu relation; indeed, this is one of the key reasons we were compelled to prove this relation in such great generality.

Surprisingly, our identification of optimizers for \eqref{1:45} makes no use of the corresponding Euler--Lagrange equation.  This constitutes a key departure from the argument employed in \cite{KV2022} in the KdV setting.  The Euler--Lagrange equation would identify $u$ as a linear combination of the \emph{square moduli of eigenfunctions}, provided that one can verify that these are linearly independent.  Instead, we combine the spectral theorem and the Wu relation to identify $u_+$ as a linear combination of the eigenfunctions themselves!  This leads us directly into the path of the following powerful spectral characterization of multisolitons:

\begin{theorem} [Spectral characterization of multisolitons, \cite{Sun2021}] \label{T: Sun}
Fix \( N\geq 1 \) and a set $\Lambda$ of negative parameters $\la_1<\cdots<\la_N$. A function $ u$ belongs to the set  $\{Q_{\Lambda, \vec c} :\, \vec c\in \R^N\}$ if and only if it satisfies the following three conditions:
\begin{SL}
\item \( \langle x \rangle u \in L^2(\mathbb{R}) \),
\item $\la_1, \ldots, \la_N$ are the only negative eigenvalues of the Lax operator \( L_u \), and 
\item \( u_+ \) belongs to the span of the corresponding \( N \) eigenfunctions.
\end{SL}
\end{theorem}

Unfortunately, condition (i) in Theorem~\ref{T: Sun} demands rather stringent spatial decay on the optimizers $u$.  Much of Section~\ref{S:VChar} is devoted to showing that optimizers enjoy such strong decay.  First, in Proposition~\ref{P:CT} we prove optimal decay of eigenfunctions associated to discrete eigenvalues.  Through the connection of optimizers to linear combinations of eigenfunctions, this guarantees some spatial decay for optimizers.  Alas, this decay is \emph{not} sufficient to apply Theorem~\ref{T: Sun}.

The solution to this impasse is a key cancellation that occurs when $u$ is reconstructed from its positive and negative frequency parts.  This discovery is elaborated in Proposition~\ref{P:1 done} and used to prove that optimizers do indeed have sufficient spatial decay to invoke Theorem~\ref{T: Sun}.

Our technique for proving spatial decay originates from the Combes--Thomas technique employed in the study of Schr\"odinger operators; see \cite{MR1785381}.  Unlike that setting, however, eigenfunctions of the \eqref{BO} Lax operator only decay polynomially, not exponentially.  The primary culprit for this is the nonlocality of the potential term originating from the Cauchy--Szeg\H{o} projection $\CS$; naturally, this also interferes with our implementation of the Combes--Thomas argument.

With a variational characterization in place, our goal is to argue in the style of Cazenave--Lions \cite{MR677997}.  In our setting, this means using a concentration compactness principle to demonstrate that every optimizing sequence $u_n$ converges to the set of optimizers, that is, to the manifold of multisolitons.

In the original setting of \cite{MR677997}, optimizing sequences are shown to converge modulo translations; subadditivity of the energy is used to show that optimizing sequences cannot split into multiple bubbles of concentration.  In our setting, such dichotomy (as Lions termed it) is actually a reality!  Indeed, as time progresses, exact multisolitons do fission into a chain of well-separated single solitons.

While it is not subadditive, we are able to prove that $\beta(\kappa;u)$ behaves additively when $u$ is comprised of well-separated bubbles and a suitable remainder; see  \eqref{eq: beta decoup} and \eqref{decoup beta}.  This is not a trivial matter.  Looking at \eqref{beta defn}, we see that this is not only a matter of showing that the resolvent $(L_u+\kappa)^{-1}$ has off-diagonal decay, but also that it depends locally on its own potential!  The need for such decoupling strongly influenced the formulation of our concentration compactness principle, Theorem~\ref{T:CC}, and proving it takes up the bulk of Section~\ref{S:CC}.

The discussion in Section~\ref{S:CC} focuses on the behavior of $\beta(\kappa;u)$ under our concentration compactness principle; we do not discuss the behavior of the eigenvalue constraints until Proposition~\ref{P:poles}.   The physical intuition is that the manifestation of each eigenvalue requires some growth of $u$, but any such growth in $u$ correspondingly increases $\beta(\kappa;u)$.  In this way, we see that every bubble of concentration must contribute solely toward satisfying the constraints, that is, creating one or more of the prescribed eigenvalues $\Lambda$.  Anything else would be wasteful.  This is what Proposition~\ref{P:poles} makes precise.  By applying it inductively, we then show that the sets of eigenvalues associated to each bubble of concentration form a partition of the set $\Lambda$.  Using our variational characterization, we are also able to demonstrate that each such bubble of concentration is a multisoliton.

Once we established everything promised in the previous paragraph, it seems that we are very close to completing the proof of Theorem~\ref{Th:  Orbital Stability}.  In fact, two important wrinkles remain.  First, we have argued that optimizing sequences must resemble a linear combination of well-separated multisolitons.  We must explain why such an object may be approximated by a single multisoliton from the manifold \eqref{MSM}.  This duty is discharged in Section~\ref{S:Multi}.  This is also the opportune place to elaborate some basic properties of multisolitons that are used elsewhere in the paper.

The second wrinkle is that for a combination of technical and expository reasons, the proof we have just outlined is all carried out for regularities $-\frac12<s<0$.  This also reflects our desire to demonstrate stability for the noisiest possible perturbations.  However, Theorem~\ref{Th:  Orbital Stability} promises more: if we begin close to a multisoliton in a higher-regularity sense, then we remain close in this \emph{same} higher-regularity sense.  This extension is carried out in Section~\ref{S:HR} via two different methods.  When $0\leq s<\frac12$, we prove that optimizing sequences are $H^s$-equicontinuous, which allows us to upgrade convergence at lower regularity to $H^s$-convergence.  When $s=\frac12$, we employ conservation of energy and argue in the spirit of the Radon--Riesz Theorem.

The overall strategy we have just outlined has definite parallels with that employed in \cite{KV2022}: both introduce a variational characterization based on a generating function for the traditional conserved quantities and eigenvalue constraints; then both employ concentration compactness techniques to show that optimizing sequences must converge to the manifold of optimizers, namely, multisolitons.

However, once we turn to the details, the analogies begin to fall apart.  Most strikingly, the Lax operators for the two hierarchies are very different.  The paper \cite{KV2022} relies on long-established tools for Sturm--Liouville operators, including the decay of eigenfunctions, locality of the resolvent, the inner/outer factorization of the reciprocal reflection coefficient, and inverse scattering theory.  These are all important for identifying optimizers and proving decoupling of profiles in the variational principle.   By contrast, $L_u$ is a fundamentally nonlocal operator whose spectral theory is very poorly understood.

In \cite{KV2022}, optimizers are seen to be Schwartz functions directly from the Euler--Lagrange equation and are then identified as multisolitons via inverse scattering.  This is not at all how we identify optimizers here; indeed, there is currently no inverse scattering theory for \eqref{BO} that can handle slowly-decaying solutions, such as solitons. Instead, in this paper, the identification of optimizers relies on establishing both the sharp decay of eigenfunctions and the Wu relation in what appears to us to be full generality (at least from the perspective of the well-posedness problem).

\section{Notation and preliminaries} \label{S:Prelims}

We adopt the following rule for the operator $\Pi$: Its precedence is lower than multiplication indicated by juxtaposition (e.g., $fg$), but higher than multiplication indicated with a dot, addition, and subtraction.  Thus,
\begin{equation*}
\CS f \cdot \CS gh + u = \bigl[\CS f\bigr]  \bigl[\Pi\bigl(gh\bigr)\bigr] + u .
\end{equation*}

Our convention for the Fourier transform is
\begin{align}\label{FT}
\widehat f(\xi) = \tfrac{1}{\sqrt{2\pi}} \int_\R e^{-i\xi x} f(x)\,dx  \quad \text{so}\quad f(x) = \tfrac{1}{\sqrt{2\pi}} \int_\R e^{i\xi x} \widehat f(\xi)\,d\xi.
\end{align}
This Fourier transform is unitary on $L^2(\R)$ and we have the Plancherel and convolution identites
\begin{align}\label{Plancherel}
    \|f\|_{L^2(\mathbb R)}=\|\widehat f\,\|_{L^2(\mathbb R)} \qtq{and}  \widehat {fg}= \tfrac{1}{\sqrt{2\pi}} \widehat f \ast \widehat g.
\end{align}

Let \(\varphi\)  be a smooth, non-negative function supported on $|\xi|\leq 2$ with $\varphi(\xi)=1$ for $|\xi|\leq 1$. For $N\in 2^\Z$ we define the Littlewood--Paley projection $P_{\leq N}$ as the Fourier multiplier operator with symbol $\varphi(\xi/N)$.  We define $P_{>N} := 1-P_{\leq N}$.  

For \(\sigma\in \R\) and \(\kappa\geq 1\) we define the Sobolev space $H^\sigma_\kappa(\R)$ as the completion of $\mathcal S(\R)$ with respect to the norm
\begin{equation}\label{H^s_k defn}
\|f\|_{H^\sigma_\kappa}^2 = \int_\R  (|\xi|+\kappa)^{2\sigma} |\widehat{f}(\xi)|^2\,d\xi .
\end{equation}
When $\kappa=1$, we simply write $H^\sigma(\R)$. Throughout the paper, we will employ the $L^2$ pairing $\langle g, f \rangle = \int \overline{ g(x)} f(x)\,dx$, which identifies $H^{\sigma}_\kappa(\R)$ and $H^{-\sigma}_\kappa(\R)$ as dual spaces.

The Hardy--Sobolev spaces are defined as 
$$
H^\sigma_+(\mathbb{R}) := \{ f \in H^\sigma(\mathbb{R}) :\, \supp(\widehat{\mkern0.3mu f}\mkern1mu) \subseteq [0,+\infty) \}.
$$

For $s>-\frac12$ we have $s+1>\frac{1}{2}$ and so the space $H^{s+1}_\kappa$ is an algebra. In particular,
\begin{equation}\label{high reg alg}
\|fg\|_{H^{s+1}_\kappa} \lesssim \|f\|_{H^{s+1}} \|g\|_{H^{s+1}_\kappa} \quad\text{uniformly for $\kappa\geq 1$.}
\end{equation}

Our next lemma contains the well-known estimate \eqref{BSm} for multipliers on lower regularity spaces. We include a proof because we employ the same ingredients to prove the spatial decoupling estimate \eqref{int}.  Recall that $W^{1,\infty}(\R)$ is comprised of bounded Lipschitz functions and is equipped with the norm
$$
\|f\|_{W^{1,\infty}}= \|f\|_{L^\infty} + \|f'\|_{L^\infty}.
$$

\begin{lemma}[Multipliers on $H^\sigma$]\label{L:multiplier}
For any $|\sigma|\leq 1,$ we have 
\begin{align}\label{BSm}
\|mf\|_{H^\sigma} \lesssim \|m\|_{W^{1,\infty}} \|f\|_{H^\sigma}.
\end{align}
Moreover, for each $\theta>\frac12$,
\begin{align}
\sup_{y\in\R} \,\bigl\|\langle x-y\rangle^{-\theta} f \bigr\|_{H^\sigma} \lesssim\|f\|_{H^\sigma}, \label{sup}\\
\int_\R\, \bigl\|\langle x-y\rangle^{-\theta} f \bigr\|_{H^\sigma}^2\,dy\simeq \|f\|_{H^\sigma}^2.  \label{int}
\end{align}
\end{lemma}

\begin{proof}
Applying the product rule, we see that
\begin{align}\label{BSm+1}
\|mf\|_{H^1} \lesssim \|m \|_{L^\infty} \|f\|_{H^1} + \|m '\|_{L^\infty}\|f\|_{L^2} \lesssim \|m\|_{W^{1,\infty}} \|f\|_{H^1} .
\end{align}
This proves \eqref{BSm} for $\sigma=1$.  The case $\sigma=-1$,
\begin{align}\label{BSm-1}
\|mf\|_{H^{-1}} \lesssim \|m\|_{W^{1,\infty}} \|f\|_{H^{-1}},
\end{align}
follows from \eqref{BSm+1} and duality.  These two estimates show that
\begin{align}\label{BSm z}
\bigl\| \langle\partial\rangle^{z} \,\bigl[m\, \langle\partial\rangle^{-z} f\bigr]\bigr\|_{L^2} \lesssim \|m\|_{W^{1,\infty}} \|f\|_{L^2}
\end{align}
uniformly for $\Re z =\pm 1$ and $\Im z\in\R$.  (It is important here that the Fourier multiplier $\langle\partial_x\rangle^{ib}$ acts unitarily on $L^2(\R)$ for any $b\in\R$.)  Complex interpolation then extends \eqref{BSm z} to the entire strip $\Re z \in [-1,1]$, which proves \eqref{BSm}.

The bound \eqref{sup} is an immediate consequence of \eqref{BSm}; indeed, this claim is valid for any $\theta\geq 0$.  By comparison, the restriction $\theta>\frac12$ is essential for \eqref{int}, even when $\sigma=0$.  In fact, the Fubini Theorem shows
\begin{align}\label{int sigma=0}
\int_\R\, \bigl\|\langle x-y\rangle^{-\theta} f \bigr\|_{L^2_x}^2\,dy  =  \int_\R \langle y\rangle^{-2\theta}\,dy \cdot  \|f\|_{L^2}^2.  
\end{align}

In order to use complex interpolation, we recast $\text{LHS}\eqref{int} \lesssim \text{RHS}\eqref{int}$ as 
\begin{align}\label{int z}
\int_\R\, \bigl\| \langle\partial\rangle^{z} \bigl[\langle x-y\rangle^{-\theta} \langle\partial\rangle^{-z} f\bigr] \bigr\|_{L^2_x}^2\,dy \lesssim  \|f\|_{L^2}^2
	\quad\text{uniformly for $\Re z \in [-1,1]$}.
\end{align}
Unitarity of imaginary powers shows that there is no extra generality here.  It also reduces our obligation to showing that \eqref{int z} is valid for $z=\pm 1$.

For any smooth function $m(x)$, the product rule yields the identities
\begin{align*}
(1+\partial) \bigl[m (1+\partial)^{-1} f\bigr] &= m f + m'  (1+\partial)^{-1} f, \\
(1+\partial)^{-1} \bigl[m (1+\partial) f\bigr] &= m f -  (1+\partial)^{-1} [m' f] .
\end{align*}
Choosing $m(x)=\langle x-y \rangle^{-\theta}$, which satisfies $|m'(x)|\lesssim m(x)$, we deduce that
\begin{align*}
\bigl\| (1+\partial) \bigl[\langle x-y \rangle^{-\theta} (1+\partial)^{-1} f \bigr]\bigr\|_{L^2}^2 &\lesssim \bigl\| \langle x-y \rangle^{-\theta} f \bigr\|_{L^2}^2 + \bigl\| \langle x-y \rangle^{-\theta} (1+\partial)^{-1}f \bigr\|_{L^2}^2,\\
\bigl\| (1+\partial)^{-1} \bigl[\langle x-y \rangle^{-\theta} (1+\partial) f \bigr]\bigr\|_{L^2}^2 &\lesssim \bigl\| \langle x-y \rangle^{-\theta} f \bigr\|_{L^2}^2.
\end{align*}
Next we integrate over $y\in\R$ and apply \eqref{int sigma=0} to obtain \eqref{int z} for $z=\pm1$.  As noted earlier, this yields the general case by complex interpolation.

It remains to verify $\text{LHS}\eqref{int} \gtrsim \text{RHS}\eqref{int}$.  Given $f\in H^\sigma(\R)$ and $g\in H^{-\sigma}(\R)$, we use Cauchy--Schwarz and that  $\text{LHS}\eqref{int} \lesssim \text{RHS}\eqref{int}$ to obtain
\begin{align*}
\biggl| \int g(x) f(x) \,dx \biggr| & \simeq \biggl| \iint  g(x) \langle x-y\rangle^{-2\theta} f(x) \,dy\,dx \biggr| \\
&\lesssim \int  \|\langle x-y\rangle^{-\theta} g\|_{H^{-\sigma}}^2  \,dy \cdot \int \|\langle x-y\rangle^{-\theta} f\|_{H^\sigma}^2\, dy  \\
&\lesssim  \|g\|_{H^{-\sigma}} \int \|\langle x-y\rangle^{-\theta} f\|_{H^\sigma}^2\, dy.
\end{align*}
Optimizing over $g\in H^{-\sigma}(\R)$ yields the desired inequality.
\end{proof}

Our next lemma recalls multiplier estimates that appeared in Lemmas~2.1 and~3.1 from \cite{Killip2024}:

\begin{lemma}[\cite{Killip2024}]\label{L:bds} For $|s|<\frac12$ and $\kappa\geq 1$ we have
\begin{align}
\bigl\|fg\|_{H^s_\kappa}& \lesssim \bigl[ \|g\|_{L^\infty} + \| g\|_{H^{1/2}} \bigr] \| f \|_{H^s_\kappa}
	\lesssim_s \kappa^{-\frac12-s} \|f\|_{H^s_\kappa}\|g\|_{H^{s+1}_\kappa},\label{ke1}\\
\bigl\|\CS u\CS f\bigr\|_{H^{-1/2}_\kappa} 
&\lesssim_s \kappa^{-\frac12-s}\| u\|_{H^s_\kappa} \|f\|_{H^{1/2}_\kappa}.\label{ke2}
\end{align}
\end{lemma}

We note that \eqref{ke2} provides the justification for condition \eqref{k0 intro} that recurs throughout this paper.  Concretely, \eqref{ke2}  shows that for each $-\frac12 < s \leq 0$, there is a constant $C_s\geq 1$ so that
\begin{align}\label{small pert}
\kappa \geq \kappa_0(u):= C_s \bigl( 1 + \|u\|_{H^s} \bigr)^{\frac2{1+2s}} \implies \| \CS uf \|_{H^{-1/2}_\kappa} \leq \tfrac12 \| f \|_{H^{1/2}_\kappa} 
\end{align}
whenever $f\in \Hohp(\R)$.  This then implies that 
\begin{align}\label{small pert'}
\tfrac32 \langle f , (L_0 +\kappa) f\rangle \geq   \langle f , (L_u +\kappa) f\rangle \geq  \tfrac12 \langle f , (L_0 +\kappa) f\rangle
	\geq \tfrac\kappa2 \| f\|_{L^2}^2 .
\end{align}
In particular, choosing $\kappa=\kappa_0(u)$, we see that 
\begin{align}\label{blue star}
 L_u  \geq - \tfrac12 C_s \bigl( 1 + \|u\|_{H^s} \bigr)^{\frac2{1+2s}}.
\end{align}

Other applications of the estimates from Lemma~\ref{L:bds} in \cite{Killip2024} include: the realization of $L_u$ as a self-adjoint operator, the identification of its form domain as $H^{1/2}_+(\R)$, and the identification of its essential spectrum.  We collect all these properties in the following result:

\begin{theorem}\label{T:Lu}
Fix $|s|<\frac12$. Given $u\in H^s(\R)$, there is a unique self-adjoint, semi-bounded operator $L_u$ with quadratic form domain $H^{1/2}_+(\R)$.  Its essential spectrum is $\sigma_\text{ess}(L_u)= [0,\infty)$. Moreover, whenever $\kappa\geq \kappa_0(u)$ with $\kappa_0(u)$ as in \eqref{small pert}, the resolvent $(L_u+\kappa)^{-1}$ exists, maps $H^{-1/2}_+(\R)$ into $H^{1/2}_+(\R)$, and satisfies
\begin{align}\label{resolvent bound}
\bigl\| (L_u+\kappa)^{-1} f\bigr\|_{H^{s+1}_\kappa} \lesssim \|f\|_{H^s_\kappa}.
\end{align}
Finally, if $v\in H^s(\R)$ and $\kappa\geq \max\{\kappa_0(u), \kappa_0(v)\}$, then 
\begin{align}\label{eq:resolvent_diff_control}
\big\|(L_u+\kappa)^{-1}-(L_v+\kappa)^{-1}\big\|_{H^s_{\kappa}\to H^{s+1}_{\kappa}}\lesssim \|u-v\|_{H^{s}_{\kappa}}.
\end{align}
\end{theorem}

\begin{proof}
All claims except \eqref{eq:resolvent_diff_control} were proved in \cite[Proposition~3.2]{Killip2024}, employing the estimates of Lemma~\ref{L:bds}.

Using the resolvent identity, \eqref{ke1}, and \eqref{resolvent bound}, we may bound
\begin{align*}
\big\| (L_u+\kappa)^{-1} - &(L_v+\kappa)^{-1} \bigr\|_{H^s_\kappa \to H^{s+1}_\kappa}\\
&= \bigl\|(L_u+\kappa)^{-1} \CS (u - v) (L_v+\kappa)^{-1} \bigr\|_{H^s_\kappa \to H^{s+1}_\kappa}\\
&\lesssim \bigl\|(L_u+\kappa)^{-1}\bigr\|_{H^s_\kappa \to H^{s+1}_\kappa} \|u - v\|_{H^s_\kappa} \bigl\|(L_v+\kappa)^{-1} \bigr\|_{H^s_\kappa \to H^{s+1}_\kappa}\\
&\lesssim \|u - v\|_{H^s_\kappa},
\end{align*}
completing the proof.
\end{proof}

Theorem~\ref{T:Lu} ensures that 
\begin{equation}\label{2:beta defn}
\beta\bigl(\kappa; u(t)\bigr) =\bigl\langle \CS u(t), \bigl(L_{u(t)}+\kappa\bigr)^{-1} \CS u(t) \bigr\rangle
\end{equation}
is well-defined for $u\in H^s(\R)$ and $\kappa\geq \kappa_0(u)$.  When $u(t)$ is a smooth solution of \eqref{BO}, it is elementary to verify that $\beta(\kappa;u)$ is conserved using \eqref{Lax eq} and \eqref{OTP}.  It is also a continuous function of $u$; indeed, \eqref{eq:resolvent_diff_control} shows that it is Lipschitz.  Thus, this conservation law extends to $H^s$-solutions via well-posedness; see  Theorem~\ref{T:KLV main}.

In truth, the conservation of $\beta$ already played a crucial role in the proof of Theorem~\ref{T:KLV main} in \cite{Killip2024}.  For example, it was employed to demonstrate $H^s$-boundedness and equicontinuity of orbits.  Recall that a bounded set $U\subseteq H^s(\R)$ is said to be \emph{equicontinuous} if
$$
\lim_{N\to\infty} \ \sup_{u\in U} \ \| P_{\geq N} u \|_{H^s} = 0, 
$$
that is, if the high-frequency components are uniformly small.

The treatment in \cite{Killip2024} focuses on the boundedness and equicontinuity of orbits when $-\frac12<s<0$; we will recapitulate the relevant results in Proposition~\ref{P:beta}.  Boundedness and equicontinuity for regularities $0\leq s<\frac12$ will be demonstrated later in Proposition~\ref{P:s pos}.

\begin{prop}\label{P:beta}
Fix $-\frac12<s<0$.  For any solution $u(t)$ to \eqref{BO} with initial data $u(0)\in H^s(\R)$, we have 
\begin{align}\label{2:APB}
	\|u(t)\|_{H^s}\lesssim \bigl(1+ \|u(0)\|_{H^s}\bigr)^{\frac{2|s|}{1+2s}}\|u(0)\|_{H^s}.
\end{align}
Moreover, if $U\subset H^s(\R)$ is bounded and equicontinuous, then so too is
$$
U^* := \{ u(t) : t\in\R \text{ and } u\text{ is a solution of \eqref{BO} with initial data $u_0\in U$}\bigr\}.
$$
Lastly, suppose $\{u_n\}_{n\geq 1}$ is $H^s$-bounded and equicontinuous and let
$$
\kappa_1 = C_s \bigl( 1 + \sup_n \|u_n \|_{H^s} \bigr)^{\frac2{1+2s}}.
$$
Then for any fixed $\kappa\geq\kappa_1$, 
\begin{equation}\label{b to 0}
  \beta(\kappa; u_n) \to 0 \ \implies \ u_n \to 0 \quad\text{in $H^s(\R)$.}
\end{equation}
\end{prop}

\begin{proof}
The bound \eqref{2:APB} and the equicontinuity claim were both proved in \cite{Killip2024}.
Indeed, they are both deduced from the equivalence
\begin{equation}\label{beta est}
\|u\|_{H^s_\kappa}^2 \lesssim_s \int_{\kappa}^\infty \varkappa^{2s} \beta(\varkappa;u)\,d\varkappa \lesssim_s \|u\|_{H^s_\kappa}^2
\end{equation}
valid for any $\kappa\geq \kappa_0(u)$.

To prove \eqref{b to 0}, we employ \eqref{small pert'}.  This shows $(L_{u_n}+\kappa) \leq \tfrac32(L_0+\kappa)$ whenever $\kappa\geq \kappa_1$.  This implies $(L_0+\kappa)^{-1} \leq \frac32 (L_{u_n}+\kappa)^{-1}$
and consequently,
\begin{equation}
\|u_+\|_{H^{-1/2}_\kappa} \leq \tfrac32 \beta(\kappa;u) .
\end{equation}
In particular, we infer that $\beta(\kappa, u_n) \to 0$ implies $u_n\to 0$ in $H^{-1/2}(\R)$.  Convergence in $H^s(\R)$ then follows from this and $H^s(\R)$-equicontinuity; see, for example, \cite[Lemma~4.2]{Killip2019}.
\end{proof}

Let us now turn our attention to the spatial decay of the resolvent $(L_u+\kappa)^{-1}$, beginning with the case $u\equiv0$.

\begin{lemma}\label{L:inverse}
For any $\kappa>0$, we have
\begin{gather}\label{E:inverse}
\bigl\| (-i\p_x + \kappa )^{-1} \bigr\|_{L^p_+(\R)\to L^p_+(\R)} \lesssim \kappa^{-1}, \quad\text{uniformly for $1\leq p\leq\infty$,}  \\
	\label{E:inverse'}
\langle x \rangle \; \bigl(\lvert \partial_x \rvert + \kappa \bigr)^{-1} \;  \langle x \rangle^{-4/3} : L^2(\R)\to L^2(\R) .
\end{gather}
Here $|\partial_x|$ denotes the Fourier multiplier operator $\widehat f(\xi) \mapsto |\xi| \widehat f(\xi).$
\end{lemma}

\begin{proof}
The common element in these two estimates is that they originate in off-diagonal decay of the convolution kernel involved. 

We first consider \eqref{E:inverse} with $\kappa=1$ and choose $m\in C^\infty(\R)$ that satisfies
\begin{equation*}
m(\xi)=(\xi+1)^{-1}\quad\text{for}\quad \xi\geq 0\qquad \text{and}\qquad  m(\xi)=0 \quad \text{for} \quad \xi\leq-\tfrac12.
\end{equation*}
When restricted to $L^p_+(\R)$, this Fourier multiplier agrees with $(-i\p_x+1)^{-1}$.

Traditional methods allow us to bound the associated convolution kernel
\begin{equation}\label{E:m}
 K(x) = \int_\R e^{i\xi x} m(\xi) \,\tfrac{d\xi}{2\pi} \qtq{as follows:}
 	|K(x) |  \lesssim \begin{cases} \log\bigl(\tfrac1{|x|}\bigr), \ & 0<|x|\leq\tfrac12\\[0.5ex] \  |x|^{-2}, \ & |x|\geq\tfrac12. \end{cases}
 \end{equation}
In the latter case, one simply integrates by parts twice.  When $x$ is small, we divide the interval of integration into two regions.  For $\xi\leq |x|^{-1}$ we simply bring absolute values inside the integral.  In the complementary region, we first integrate by parts once, then bring absolute values inside the integral.

For general $\kappa>0$, the operator $(-i\p_x+\kappa)^{-1}$ has convolution kernel $K( \kappa x)$.  Using \eqref{E:m}, we obtain an $L^1$ bound on this kernel that proves \eqref{E:inverse}.

Turning to \eqref{E:inverse'} and following the same argument as above, we find that the convolution kernel 
\begin{align*}
K(x) = \int_\R \tfrac{e^{i\xi x}}{|\xi|+\kappa}\, \tfrac{d\xi}{2\pi} = \int_\R \tfrac{e^{i\xi \kappa x}}{|\xi|+1}\, \tfrac{d\xi}{2\pi}\qtq{satisfies} 
|K(x) |  \lesssim \begin{cases} \log\bigl(\tfrac1{|\kappa x|}\bigr), \ & 0<|\kappa x|\leq\tfrac12\\[0.5ex] \ |\kappa x|^{-2}, \ & |\kappa x|\geq\tfrac12. \end{cases}
\end{align*}
To treat the full operator appearing in \eqref{E:inverse'}, we first verify that for fixed $\kappa>0$,
\begin{align*}
A&:= \sup_x \int \langle x \rangle^{4/3} \bigl| K(x-y)\bigr| \langle y \rangle^{-4/3} \,dy < \infty, \\
B&:= \sup_y \int \langle x \rangle^{2/3} \bigl| K(x-y)\bigr| \langle y \rangle^{-4/3} \,dx < \infty.
\end{align*}
Combining this with Cauchy--Schwarz, we get
\begin{align*}
\bigl| \bigl\langle f,  \langle x \rangle \bigl(\lvert \partial_x \rvert + \kappa \bigr)^{-1} \langle x \rangle^{-4/3}  g\bigr\rangle \bigr|
	&\lesssim \iint |f(x)| \langle x \rangle  \bigl| K(x-y)\bigr| \langle y \rangle^{-4/3} |g(y)| \,dy\,dx \\
	&\lesssim \sqrt{AB} \, \| f\|_{L^2} \| g\|_{L^2}.
\end{align*}
This completes the proof of \eqref{E:inverse'}.
\end{proof}

It is physically intuitive that the total energy of two widely separated wavepackets is well approximated by the sum of the two energies; indeed, similar reasoning applies also to mass and momentum.  In Section~\ref{S:CC}, we will need this type of decoupling for  $\beta(\kappa; u)$.  This is rather more subtle.  One part of this relies on off-diagonal decay of the resolvent $(L_u+\kappa)^{-1}$.  This is the topic of Proposition~\ref{P:locality}.   We will be employing a Combes--Thomas style of argument, which is a common technique for differential operators.  However, to overcome the nonlocality of our Lax operator, we rely on the following result from the Calderon--Zygmund theory:

\begin{lemma}\label{L:AP}
Fix $|s|\leq 1$ and $1<p<\infty$.  We have
\begin{align}
\bigl\|(1-i\rho x)^\theta \CS (1-i\rho x)^{-\theta} \bigr\|_{L^p\to L^p} &\lesssim_{p,\theta} 1 \qquad\text{for}\quad -\tfrac1p<\theta<1-\tfrac1p,\label{H Lp}\\
\bigl\| (1-i\rho x)^\theta \CS (1-i\rho x)^{-\theta}\bigr\|_{H^s\to H^s}&\lesssim_{\theta} 1 \qquad\text{for} \quad -\tfrac12<\theta<\tfrac12,\label{H Hs}
\end{align}
uniformly for $0\leq \rho\leq 1$.
\end{lemma}

\begin{proof}
For $-\frac1p<\theta<1-\tfrac1p$ and $\rho\in \R$, $|1-i\rho x|^{p\theta}$ is an $A_p$ weight; see \cite[Ch.~V]{bigStein}.  This implies that Calder\'on--Zygmund operators (such as the Hilbert transform) are bounded on $L^p( \R;|1-i\rho x|^{p\theta}\,dx)$.  As a consequence, $\Pi=\frac12(\Id+iH)$ satisfies \eqref{H Lp}.  A rescaling argument shows that the implicit constant is independent of $\rho$. 

Turning now to \eqref{H Hs}, we have 
\begin{align*}
\partial_x \bigl[(1-i\rho x)^{\theta} \CS  (1-i\rho x)^{-\theta}\bigr]&=  (1-i\rho x)^{\theta} \CS (1-i\rho x)^{-\theta}\partial_x \\
&\quad - \tfrac{i\rho\theta}{1-i\rho x} (1-i\rho x)^{\theta} \CS  (1-i\rho x)^{-\theta}\\
&\quad +(1-i\rho x)^{\theta} \CS  (1-i\rho x)^{-\theta}\tfrac{i\rho\theta}{1-i\rho x}.
\end{align*}
Applying \eqref{H Lp} with $p=2$ and $|\theta|<\frac12$, we deduce that
\begin{align}\label{weighted bdd on H1}
\bigl\|(1-i\rho x)^{\theta} \CS  (1-i\rho x)^{-\theta}\bigr\|_{H^1\to H^1} \lesssim_\theta 1,
\end{align}
uniformly for $0\leq \rho\leq 1$.  A simple duality argument then gives
\begin{align}\label{weighted bdd on H-1}
\bigl\|(1-i\rho x)^{\theta} \CS  (1-i\rho x)^{-\theta}\bigr\|_{H^{-1}\to H^{-1}} \lesssim_\theta 1.
\end{align}

Claim \eqref{H Hs} follows from \eqref{weighted bdd on H1} and \eqref{weighted bdd on H-1} via complex interpolation.  Indeed, 
$$
\langle \partial\rangle ^{z} (1-i\rho x)^{\theta} \CS  (1-i\rho x)^{-\theta}\langle \partial\rangle ^{-z}
$$
is a holomorphic family of operators on the strip $\Re z \in [-1,1]$ that is $L^2$-bounded on the boundary of the strip in view of \eqref{weighted bdd on H1} and \eqref{weighted bdd on H-1}.  It is important here that $\langle \partial\rangle ^{ib}$ is unitary on $L^2(\R)$ for all $b\in \R$.
\end{proof}

\begin{prop}[Locality of the resolvent]\label{P:locality}
Fix $|s|<\frac12$, $|\theta|<\frac12$, and $M>0$.  There exists a constant $C_{s,\theta}\geq 1$ so that
$$
\bigl\| (1-i\rho x)^\theta (L_u+\kappa)^{-1}(1-i\rho x)^{-\theta}\|_{H^s_+\to H^{s+1}_+} \lesssim 1,
$$
uniformly for $0\leq \rho\leq1$, $u\in H^s(\R)$ with $\|u\|_{H^s}\leq M$, and $\kappa\geq C_{s, \theta} ( 1 + M)^{\frac2{1+2s}}$.
\end{prop}

\begin{proof}
By duality, it suffices to consider $0\leq\theta<\frac12$. Fix $u\in H^s(\R)$ with $\|u\|_{H^s}\leq M$.  For $\kappa\geq 1$ sufficiently large, we will construct 
\begin{align*}
(1-i\rho x)^\theta (L_u+\kappa)^{-1}(1-i\rho x)^{-\theta} &=\Bigl[(1-i\rho x)^\theta (L_u+\kappa)(1-i\rho x)^{-\theta}\Bigr]^{-1}\\
&= \Bigl[-i\partial + \kappa + \tfrac{\rho\theta}{1-i\rho x} \!-\! (1-i\rho x)^\theta \CS u (1-i\rho x)^{-\theta}\Bigr]^{-1}
\end{align*}
as a bounded operator from $H^s_+(\R)$ to $H^{s+1}_+(\R)$ via a Neumann series, employing the observation that 
$$
\bigl\|(-i\partial + \kappa)^{-1}\bigr\|_{H^s_+\to H^{s+1}_+} \lesssim 1.
$$
For this purpose, it suffices to show that we may choose $C_{s,\theta}$ sufficiently large so that 
\begin{align}\label{goal 101}
\bigl\| \bigl[\tfrac{\rho\theta}{1-i\rho x} \!-\!(1- i\rho x)^\theta \CS \,u (1-i\rho x)^{-\theta}\bigr](-i\partial + \kappa)^{-1}\bigr\|_{H^s_+\to H^s_+} <1
\end{align} 
for all $\kappa\geq C_{s, \theta} ( 1 + M)^{\frac2{1+2s}}$ and $0\leq\rho\leq1$.

Using \eqref{BSm}, for $0\leq \rho\leq1$ we have
\begin{align*}
\bigl\| \tfrac{\rho\theta}{1-i\rho x} (-i\partial + \kappa)^{-1}\bigr\|_{H^s_+\to H^s_+}
\lesssim \bigl\|\tfrac{\rho\theta}{1-i\rho x}\bigr\|_{W^{1,\infty}} \bigl\| (-i\partial + \kappa)^{-1}\bigr\|_{H^s_+\to H^s_+}\lesssim \kappa^{-1},
\end{align*} 
which is acceptable, provided $C_{s,\theta}$ is chosen sufficiently large.

For $f\in H^s(\R)$ and $0\leq \rho\leq 1$, we use \eqref{H Hs} and  \eqref{ke1} to bound
\begin{align*}
\bigl\| (1-i\rho x)^\theta& \CS \,u (1-i\rho x)^{-\theta}(-i\partial + \kappa)^{-1}f\bigr\|_{H^s}\\
&\lesssim \bigl\| (1-i\rho x)^\theta \Pi(1-i\rho x)^{-\theta}\bigr\|_{H^s\to H^s} \bigl\|u(-i\partial + \kappa)^{-1}f\bigr\|_{H^s}\\
&\lesssim_{s, \theta} \|u\|_{H^s}\bigl[ \|(-i\partial + \kappa)^{-1}f\|_{L^\infty} + \|(-i\partial + \kappa)^{-1}f\|_{H^{\frac12}}\bigr]\\
&\lesssim_{s, \theta} M \kappa^{-\frac12 -s}\|f\|_{H^s},
\end{align*}
which is acceptable for $\kappa\geq C_{s, \theta} ( 1 + M)^{\frac2{1+2s}}$, provided $C_{s,\theta}$ is chosen sufficiently large.
\end{proof}

\begin{lemma} \label{L:Riesz} Fix $-\frac12<\sigma<s<0$ and $\rho\in\R$.  Then
$f \mapsto (1-i\rho x)^{-1} f$ is a compact operator from $H^s(\R)$ to $H^\sigma(\R)$.
\end{lemma}

\begin{proof}
For any $M\in 2^\N$, the Hilbert--Schmidt test shows that $f \mapsto (1-i\rho x)^{-1} P_{\leq M}f$ is a compact operator from $H^s(\R)$ to $H^\sigma(\R)$.  On the other hand,
\begin{align*}
\bigl\|(1-i\rho x)^{-1} P_{>M} f\bigr\|_{H^\sigma} \lesssim \bigl\|(1-i\rho x)^{-1} \bigr\|_{W^{1,\infty}} \|P_{>M}f\|_{H^\sigma} \lesssim_\rho  M^{\sigma-s} \|f\|_{H^s}.
\end{align*} 
Sending $M\to \infty$ we see that $f \mapsto (1-i\rho x)^{-1} f$ is a compact operator because it is the norm-limit of compact operators.
\end{proof}


\section{Multisoliton solutions}\label{S:Multi}

In order to make sense of Matsuno's formula \eqref{E:Qbc} for multisolitons, one must ensure that the matrix $M$ defined in \eqref{matrix A} satisfies $\det M \neq 0$.  This is well-known and will be verified as part of our next lemma.

The spatial derivative in \eqref{E:Qbc} obviates any ambiguity arising from different choices of the branch of the complex logarithm.  Indeed, taking this derivative allows us to further simplify the formula:
\begin{align}\label{Matsuno multi}
Q_{\Lambda, \vec c}\, (x)= - 2\Im \tfrac{d}{dx} \ln \det\bigl[M(x)\bigr] 
&= - 2\Im \tr \bigl\{M'(x) M(x)^{-1}\bigr\}\notag\\
&=2\Re  \tr \bigl\{M(x)^{-1}\bigr\}.
\end{align}

\begin{lemma}\label{L:M^-1}
The matrix $M(x)$ defined in \eqref{matrix A} is invertible.  Moreover,
\begin{gather}
\|M(x)^{-1}\| \lesssim \max_{1\leq j\leq N} \tfrac{1}{\langle x-c_j\rangle}, \label{M^-1}\\
\bigl| \partial_x^m Q_{\Lambda, \vec c}\, (x) \bigr| \lesssim \max_{1\leq j\leq N} \tfrac{1}{\langle x-c_j\rangle^{m+2}} \label{Q decay}
\end{gather}
for each integer $m\geq 0$.  The implicit constants here are uniform in $\vec c$, but do depend nontrivially on $m$ and the set $\Lambda$.
\end{lemma}

\begin{proof}
By the structure of $M(x)$, for any vector $\vec z\in \C^N$ we have
$$
\Re\langle \vec{z}, M(x) \vec{z}\,\rangle  = \sum_{j=1}^N \tfrac{1}{2|\la_j|} |z_j|^2\geq \min_{1\leq j\leq N} \tfrac{1}{2|\la_j|}\|\vec z{\mkern 2mu}\|^2.
$$
This shows that $M(x)$ is invertible and
\begin{align}\label{easy M^-1}
\|M(x)^{-1}\|\leq \max_{1\leq j\leq N} 2|\la_j|, \quad \text{uniformly for $x\in \R$}.
\end{align}

To continue, we decompose $M(x)$ into its diagonal part $D(x)$ and its off-diagonal part $E$. Evidently,
\begin{align*}
\|D(x)^{-1}\|\leq \max_{1\leq j\leq N} \tfrac1{|x-c_j|},
\end{align*}
while by the Hilbert--Schmidt test,
\begin{align*}
\|E\|^2\leq \sum_{j\neq k} \tfrac{1}{|\la_j-\la_k|^2} .
\end{align*}

When $ |x- c_j| > \|E\|+1$ for each $1\leq j\leq N$, we have $\|E D(x)^{-1}\|<1$ and so we can invert $M(x)$ via a Neumann series based on
$$
M(x)^{-1} = D(x)^{-1} ( 1+ ED(x)^{-1})^{-1}.
$$
Form this, we deduce
$$
\|M(x)^{-1}\|\leq 2\max_{1\leq j\leq N} \tfrac1{|x-c_j|} \simeq \max_{1\leq j\leq N} \tfrac1{\langle x-c_j\rangle}.
$$
This proves \eqref{M^-1} when $x$ is far from the centers $c_j$.  In the complementary region, \eqref{M^-1} follows directly from \eqref{easy M^-1}.

Turning to \eqref{Q decay}, we rewrite \eqref{Matsuno multi} as
\begin{align}\label{Q re}
Q_{\Lambda, \vec c}\, (x) = \tr \bigl\{M(x)^{-1} + M(x)^{-\dagger} \bigr\} = \tr \bigl\{M(x)^{-1}\bigl[M(x)+ M(x)^{\dagger}\bigr]M(x)^{-\dagger} \bigr\},
\end{align}
where $\dagger$ is used to denote the  conjugate transpose.  Notice that $M+ M^{\dagger}$ is the diagonal matrix with entries $|\lambda_j|^{-1}$; it is bounded and it is independent of $x$.  In this way, the $m=0$ case of \eqref{Q decay} follows directly from \eqref{M^-1}.

We now differentiate \eqref{Q re}.  From \eqref{matrix A}, we see that $\partial_x M^{-1} = i M^{-2}$ and so
\begin{align*}
\bigl| \partial_x^m Q_{\Lambda, \vec c}\, (x) \bigr| \lesssim_m \bigl\| M(x)^{-1} \bigr\|^{m+2} .
\end{align*}
Thus, \eqref{Q decay} follows from \eqref{M^-1}.
\end{proof}

The main result in this section is the following molecular structure of multisolitons, which asserts that linear combinations of well-separated multisolitons lie close to the multisoliton manifold.  We will use $\coprod$ to represent concatenation of vectors.

\begin{prop}[Molecular structure of multisolitons]\label{prop:MD} Fix $J\geq 1$.  For each $1\leq j\leq J$, fix $N_j\geq 1$, vectors $\vec{c}_j\in\R^{N_j}$, and sets $\Lambda_j\subset (-\infty, 0)$ of $N_j$ elements.  We assume that the elements of all $\Lambda_j$ are pairwise distinct. Let $N=\sum_{j=1}^JN_j$ and $\Lambda=\cup\Lambda_j.$ For each $J$-tuple of sequences $(x_n^j)_{j=1}^{J}$ satisfying 
\begin{align}\label{5:26}
\lim_{n \to \infty} (x_n^j - x_n^k) = \infty \qquad \text{whenever } j\neq k,
\end{align}
there exists a sequence of shift parameters $\vec{c}_n = \coprod (\vec{c}_{j} + x_n^j \vec{\mathbf{1}})\) such that
\be\label{molecular decomposition}
Q_{\Lambda, \vec{c}_n}(x) - \sum_{j=1}^{J} Q_{\Lambda_j, \vec{c}_j}(x - x_n^j) \longrightarrow 0 \quad \text{in $H^\sigma(\mathbb{R})$}
\ee
as $n\to\infty$ for any regularity $\sigma\in\R$.
\end{prop}

\begin{proof}
Let \( M_n:= M_{\Lambda, \vec{c}_n} \) and \( M_{\Lambda_j, \vec{c}_j} \) be as defined in~\eqref{matrix A}, so that
\begin{equation*}
        Q_{\Lambda, \vec{c}_n}(x) = 2\Re  \tr \bigl\{ M_n(x)^{-1} \bigr\}
        \ \ \text{and} \ \ %
        Q_{\Lambda_j, \vec{c}_j}(x- x_n^j) = 2\Re  \tr \bigl\{  M_{\Lambda_j, \vec{c}_j}(x- x_n^j)^{-1} \bigr\}.
\end{equation*}
Using Lemma~\ref{L:M^-1}, we have
\begin{equation}\label{M size}
\bigl\| M_n(x)^{-1} \bigr\| \lesssim 1 \qtq{and} \bigl\|  M_{\Lambda_j, \vec{c}_j}(x- x_n^j)^{-1} \bigr\| \lesssim  \langle x- x_n^j \rangle^{-1} .
\end{equation}
The implicit constants here depend on both $\Lambda$ and $\vec c_j$. What is important is that they are uniform in $x$ and $n$.

Next, we introduce the block-diagonal matrix $\widetilde{M}_n(x) \in \mathbb{C}^{N \times N}$, whose diagonal blocks are given by $M_{\Lambda_j, \vec{c}_j}(x - x_n^j)$ so that
\begin{equation*}
        \sum_{j=1}^{J} Q_{\Lambda_j, \vec{c}_j}(x - x_n^j) = 2\Re  \tr \Bigl\{  \widetilde M_n(x)^{-1} \Bigr\}.
\end{equation*}
In this way, claim \eqref{molecular decomposition} reduces to showing that
\begin{equation}\label{Re tr}
\Re \tr\bigl\{ M_n(x)^{-1} - \widetilde{M}_n(x)^{-1}\bigr\} \longrightarrow 0 \quad \text{in $H^\sigma(\mathbb{R})$} \quad\text{as $n \to \infty$.}
\end{equation}
In pursuit of this goal, we use the resolvent identity to decompose
\begin{align}\label{6.4}
\widetilde{M}_n^{-1} - M_n^{-1}
&= \widetilde{M}_n^{-1}  (M_n - \widetilde{{M}}_n) M_n^{-1}\notag\\
&=\widetilde{M}_n^{-1}  (M_n - \widetilde{{M}}_n) \widetilde{M}_n^{-1} \bigl[ 1 - (M_n - \widetilde{M}_n)M_n^{-1}\bigr].
\end{align}

The matrix $\widetilde M_n(x)$ agrees with $M_n(x)$ on each of the $N_j\times N_j$ diagonal blocks.  This has a number of important consequences.  First, the matrix
$M_n(x) - \widetilde{{M}}_n(x)$ is independent of both $n$ and $x$.  Correspondingly, it is easy to see that
\begin{equation}\label{sup diff}
 \sup_{n,x} \  \bigl\| M_n(x) - \widetilde{{M}}_n(x) \bigr\| \lesssim 1.
\end{equation}
Secondly, by \eqref{M size},
\begin{align}\label{first term}
\bigl\|\widetilde{M}_n(x)^{-1}[M_n(x) - \widetilde{M}_n(x)]\widetilde{M}_n(x)^{-1}\bigr\|
\lesssim \sum_{j\neq k}\frac{1}{\langle x-x_n^k\rangle\langle x-x_n^j\rangle}.
\end{align} 
Combining \eqref{M size} and \eqref{sup diff}, we obtain
\begin{align}\label{5:21}
\bigl\|  1 - [M_n(x) - \widetilde{M}_n(x)]M_n(x)^{-1}\bigr\| \lesssim 1,
\end{align}
uniformly in $x$ and $n$.

Putting together \eqref{6.4}, \eqref{first term}, and \eqref{5:21}, and employing \eqref{5:26}, we deduce that
\begin{equation}\label{L2 zero}
\int_\R\bigl\|  \widetilde{M}_n(x)^{-1} - M_n(x)^{-1} \bigr\|^2\, dx \lesssim  \int_\R\sum_{j\neq k}\frac{1}{\langle x-x_n^k\rangle^2\langle x-x_n^j\rangle^2}\, dx
	\longrightarrow 0
\end{equation}
as $n\to\infty$.  This proves \eqref{Re tr} with $\sigma=0$, which then covers all $\sigma\leq 0$.  Similarly, we can cover general $\sigma$ by proving \eqref{Re tr} only for $\sigma\in\N$.

As noted in the proof of Lemma~\ref{L:M^-1}, $\partial_x M^{-1} = i M^{-2}$.  Thus, for $\sigma\in\N$,
\begin{align*}
\bigl\|  \partial_x^\sigma \bigl[ \widetilde{M}_n(x)^{-1} - M_n(x)^{-1}\bigr] \bigr\|
	&\lesssim_\sigma \bigl\| \widetilde{M}_n(x)^{-\sigma-1} - M_n(x)^{-\sigma-1} \bigr\| \\
&\lesssim_\sigma \bigl[ \| \widetilde{M}_n^{-1}(x) \| + \| M_n^{-1}(x) \|\bigr]^{\sigma}\bigl\| \widetilde{M}_n(x)^{-1} - M_n(x)^{-1} \bigr\|.
\end{align*}
Claim \eqref{Re tr} for $\sigma\in\N$ then follows from \eqref{L2 zero} and \eqref{M size}.
\end{proof}

As we will employ Sun's spectral characterization of multisoliton profiles Theorem~\ref{T: Sun}, we must address the fact the exact formula used in \cite{Sun2021} differs meaningfully from \eqref{Matsuno multi}. Accounting for differences in notation, \cite[Corollary~1.10]{Sun2021} describes the multisoliton profiles as 
\begin{align}\label{Sun multi}
Q_{\Lambda, \vec c}\, (x)=2\Re  \sum_{1\leq j, k\leq N} \bigl[M(x)^{-1}\bigr]_{jk}.
\end{align}
The equivalence of the representations \eqref{Matsuno multi} and \eqref{Sun multi} follows from the following observation: If $A = \diag(\la_1, \ldots, \la_N)$, then the commutator $[A, M]$ has entries
$$
[A,M]_{jk} = \delta_{jk}-1.
$$
Cycling the trace, this yields
$$
\tr\bigl\{M^{-1}\bigr\} - \sum_{1\leq j, k\leq N} \bigl[M(x)^{-1}\bigr]_{jk} = \tr\bigl\{M^{-1}[A,M]\bigr\}= 0.
$$

\section{The Wu relations}\label{S:Wu}

In his original work, Wu proved the key relation \eqref{Wu relation Hs} for $\lambda<0$ with quite strong hypotheses on $u$; see \cite[Lemma~2.5]{MR3484397}.  His approach was by direct computation in Fourier variables.  Wu's hypotheses were subsequently relaxed in Sun \cite[Proposition~2.4]{Sun2021}, who also introduced an operator-theoretic perspective by connecting the identity \eqref{Wu relation Hs} to a commutator relation between the Lax operator $L_u$ and the operator $X$ defined thus:

\begin{defi}\label{D:X}
Let $X$ denote the (unbounded) operator on $L^2_+(\R)$ defined as the generator of the contraction semigroup
\begin{equation}\label{E:D:X}
e^{-itX} f = \tfrac{1}{\sqrt{2\pi}} \int_0^\infty e^{i\xi x} \widehat f(\xi + t)\,d\xi = \Pi\bigl( e^{-itx} f \bigr)
\end{equation}
acting on $L^2_+(\R)$.  Concretely,
$$
D(X) = \bigl\{ f\in L^2_+ (\R): \widehat f\in H^1\bigl([0,\infty)\bigr) \bigr\} \quad \text{and} \quad \widehat{Xf}(\xi) =i\tfrac{d{\widehat f}}{d\xi}(\xi)
	\quad \text{for} \quad f\in D(X).
$$
\end{defi}

If $f\in L^2_+(\R)$ satisfies $\langle x\rangle f(x)\in L^2(\R)$, then $f\in D(X)$ and $[Xf](x)=xf(x)$.  We have elected to call the operator $X$ in deference to his relation.  Although $D(X)$ is strictly larger than the domain of the operator $f\mapsto x f$, it still presents a major obstruction to applying Sun's method to our rough and slowly decaying potentials.  Instead, we will base our approach around the semigroup \eqref{E:D:X}.  This has much better mapping properties; for example, we have the following:

\begin{lemma}\label{semi Hs}
For any $\sigma\in\R$ and $t\geq 0$,
\begin{align}\label{e^itX duh}
\bigl\| e^{-itX} f \bigr\|_{H^\sigma(\R)} \lesssim \langle t \rangle^{|\sigma|} \| f\|_{H^\sigma(\R)} .
\end{align}
\end{lemma}

\begin{proof}  This is very elementary:
\begin{align*}
\hspace*{-1.3em}
\bigl\| e^{-itX} f \bigr\|_{H^\sigma}^2
	&= \int_0^\infty \langle \xi\rangle^{2\sigma} \bigl| \widehat f(\xi+t) \bigr|^2\,d\xi
	\leq 2^{2|\sigma|}\int_0^\infty \langle t\rangle^{2|\sigma|} \langle \xi+ t\rangle^{2\sigma} \bigl| \widehat f(\xi+t) \bigr|^2\,d\xi.
 \qedhere
\end{align*}
\end{proof}

The Sobolev space $\Hohp(\R)$ will be especially important because it is the quadratic form domain of the Lax operator $L_u$.  In particular, all eigenfunctions belong to this space.  In our application, the function $f$ in our next lemma will be chosen to be such an eigenfunction and $u$ will be real-valued.

\begin{lemma}\label{L:fu}
For $u\in H^s(\R)$ with $|s|<\frac12$ and $f\in H^{1/2}(\R)$, the map $\xi\to \widehat{uf}(\xi)$ is continuous throughout the real line; moreover,
\begin{align}\label{fu bnd}
\bigl| \widehat{uf}(\xi) \bigr| \lesssim \langle\xi\rangle^{-s}  \|u\|_{H^s} \|f\|_{H^{1/2}}   \qtq{and} \sqrt{2\pi}\,\widehat{uf}(0) = \langle \overline u, f\rangle.
\end{align}
\end{lemma}

\begin{proof}
By Cauchy--Schwarz
\begin{equation}\label{1/2 and s}
\begin{aligned}
\int \bigl| \widehat u(\xi-\eta)  \bigr| \bigl|\widehat f (\eta)\bigr| \,d\eta
	&\lesssim  \sup_{\eta\in\R}\  \langle\xi-\eta\rangle^{-s}  \langle\eta\rangle^{-\frac12} \ \|u\|_{H^s} \|f\|_{H^{1/2}} \\
&\lesssim  \langle\xi\rangle^{-s}  \|u\|_{H^s} \|f\|_{H^{1/2}} .
\end{aligned}
\end{equation}
This proves the first part of \eqref{fu bnd}.

If $u$ and $f$ are Schwartz functions, then $\widehat{uf}$ is clearly continuous.  The estimate \eqref{1/2 and s} allows us to extend this to general $u$ and $f$ by approximation.

Likewise, the second relation in \eqref{fu bnd} is evident for Schwartz functions.  This extends to general $u$ and $f$ by the joint continuity of the duality relation between $H^{\pm1/2}(\R)$ and the embedding $H^{s}(\R)\hookrightarrow H^{-1/2}(\R)$.
\end{proof}

\begin{proof}[Proof of Theorem~\ref{th: Wu relation}]
Setting $L_0=-i\partial_x$ and recalling \eqref{ke2}, we know that both
\begin{align}\label{s.a. bound}
f\mapsto L_0 f \ \ \text{and}\ \ f \mapsto \CS (u f) \ \ \text{are bounded operators} \ \ H^{1/2}_+(\R) \to H^{-1/2}_+(\R).
\end{align}

By Lemma~\ref{semi Hs}, we also know that $e^{-itX}$ is bounded on both $H^{\smash[b]{\pm1/2}}_+(\R)$.  Thus, we may define
\begin{align}
\Phi_0(t) &:=  \tfrac{1}{t} \langle f, e^{-itX} L_0 f-  L_0 e^{-itX} f \rangle, \\
\Phi_u(t) &:=  \tfrac{1}{t} \langle f, e^{-itX} \Pi\, uf -  \Pi\, u\, e^{-itX}f \rangle 
\end{align}
for any $t>0$ and any $f\in \Hohp(\R)$. When $f$ is an eigenfunction of the (self-adjoint!) operator $L_u$, these combine nicely: 
\begin{align}\label{Phi=0}
\Phi_0(t) - \Phi_u(t) = \tfrac{1}{t} \langle f, e^{-itX} L_u f -  L_u e^{-itX} f \rangle = 0 \qtq{for all} t>0.
\end{align}
We will prove the Wu relation \eqref{Wu relation Hs} by considering the $t\to 0$ limit of each term on LHS\eqref{Phi=0}.

Working in Fourier variables, we have
\begin{align*}
\Phi_0(t) = \tfrac{1}{t} \int_0^\infty \overline{\widehat f (\xi)} \bigl[ (\xi+t) - \xi \bigr]  \widehat f (\xi+t) \,d\xi
	= \int_0^\infty \overline{\widehat f (\xi)} \widehat f (\xi+t) \,d\xi
\end{align*}
for any $f\in \Hohp(\R)$.  The continuity of translations in $L^2(\R)$ then shows that
\begin{align}\label{7:30}
2\pi\Phi_0(t)  \longrightarrow 2\pi\|f\|_{L^2}^2 \qtq{as} t\searrow 0.
\end{align}

For any $f\in \Hohp(\R)$ and $\xi\geq 0$, we have
\begin{align*}
\bigl[e^{-itX} \Pi\, uf -  \Pi\, u\, &e^{-itX}f\bigr]\,\widehat{\ \vphantom{M}} (\xi) \\
&= \tfrac1{\sqrt{2\pi}} \int_0^\infty \widehat u(\xi+t-\eta)\widehat f (\eta) \,d\eta - \tfrac1{\sqrt{2\pi}}\int_0^\infty \widehat u(\xi-\eta)\widehat f (\eta+t) \,d\eta \\
&= \tfrac1{\sqrt{2\pi}} \int_0^\infty \widehat u(\xi+t-\eta)\widehat f (\eta) \,d\eta - \tfrac1{\sqrt{2\pi}}\int_t^\infty \widehat u(\xi+t-\eta)\widehat f (\eta) \,d\eta \\
&= \tfrac1{\sqrt{2\pi}} \int_0^t \overline{\widehat u(\eta-t-\xi)} \, \widehat f (\eta) \,d\eta.
\end{align*}
This relied solely on a change of variables and the fact that $u$ is real-valued.  The absolute convergence of these integrals is guaranteed by \eqref{1/2 and s}.

Combining this identity and \eqref{Plancherel}, we obtain
\begin{align}\label{203}
\Phi_u(t)& = \tfrac{1}{\sqrt{2\pi}\,t}\! \int_0^\infty \!\int_0^t \; \overline{\widehat f (\xi) \, \widehat u(\eta-t-\xi)} \, \widehat f (\eta) \,d\eta \, d\xi\notag\\
	&= \tfrac{1}{t} \int_0^t \; \overline{ \widehat{ uf} (\eta-t) }  \, \widehat f (\eta) \,d\eta.
\end{align}
Absolute convergence of the integrals is again guaranteed by \eqref{1/2 and s}.

As $f$ is an eigenfunction of $L_u$ with eigenvalue $\lambda$,
\begin{align}\label{hat eigenfn}
\lambda \widehat f(\xi) = \xi \widehat f(\xi) - \widehat{ uf} (\xi) \qtq{for a.e.} \xi\geq 0.
\end{align}
Using \eqref{hat eigenfn} in \eqref{203}, we deduce that
\begin{align}\label{202}
\lambda \Phi_u(t) &= \tfrac{1}{t} \int_0^t \, \overline{ \widehat{ uf} (\eta-t) }  \, \eta \widehat f (\eta) \,d\eta  
	- \tfrac{1}{t} \int_0^t \, \overline{ \widehat{ uf} (\eta-t) }  \, \widehat{ uf} (\eta) \,d\eta.
\end{align}
Using the results of Lemma~\ref{L:fu}, we then easily deduce that
\begin{align}\label{7:31}
2\pi\lambda \Phi_u(t) \longrightarrow -  2\pi\bigl| \widehat{ uf} (0)\bigr|^2 = - \bigl| \langle u, f\rangle\bigr|^2 \qtq{as} t\searrow 0.
\end{align}

The Wu relation \eqref{Wu relation Hs} follows by combining \eqref{Phi=0}, \eqref{7:30}, and \eqref{7:31}.

Positive eigenvalues are clearly inconsistent with \eqref{Wu relation Hs}. Next, we will preclude the possibility that $\lambda=0$ is an eigenvalue. Indeed, using \eqref{hat eigenfn} with $\lambda =0$ in \eqref{203}, we find the following replacement for \eqref{202}:
\begin{align*}
\Phi_u(t) &= \tfrac{1}{t} \int_0^t \, (\eta-t) \overline{ \widehat{ f} (\eta-t) }   \widehat f (\eta) \,d\eta. 
\end{align*}
This clearly converges to zero as $t\searrow 0$.  Combining this with \eqref{Phi=0} and \eqref{7:30}, we conclude that $\|f\|_{L^2}=0$, a contradiction.

Thus, all eigenvalues are necessarily negative. To see that they are also simple, we note that if an eigenvalue were to have multiplicity $\geq 2$, then one would be able to find a non-zero eigenvector in the kernel of the linear functional $f\mapsto \langle u,f \rangle$, thereby contradicting \eqref{Wu relation Hs}.

It remains to prove \eqref{Wu 2} and that $\widehat f$ is continuous on $[0,\infty)$. As eigenfunctions belong to the form domain of the operator, we have $f\in H^{1/2}_+(\R)$. By Lemma~\ref{L:fu}, we then obtain that $\widehat{uf}(\xi)$ is continuous on $[0,\infty)$. Recasting \eqref{hat eigenfn} as
\begin{equation*}
\widehat f(\xi) =  \tfrac{1}{\xi - \lambda}  \widehat{ uf} (\xi) 1_{[0,\infty)}(\xi) 
\end{equation*}
and recalling that $\lambda<0$, we conclude that $\widehat f$ is continuous on $[0,\infty)$.  Finally, sending $\xi\searrow 0$ yields \eqref{Wu 2}.
\end{proof}

\section{Variational characterization of multisolitons}\label{S:VChar}
The goal of this section is to prove Theorem~\ref{Th: variational problem}.  Correspondingly, we consider only $-\frac12<s<0$; indeed, the result is strongest when $s$ is smallest, so there is no advantage to allowing $s\ge 0$. Throughout we fix $N\geq1$ and spectral parameters $\la_1<\cdots<\la_N<0$. 

As discussed in Section~\ref{S:Prelims}, if $u\in H^s(\R)$ then $L_u$ is self-adjoint,
\begin{equation}\label{4.1}
\sigma_{\text{ess}}(L_u) = [0, \infty), \qtq{and}    L_u\geq -\tfrac12 C_s(1+\|u\|_{H^s})^{\frac2{1+2s}}.
\end{equation}
Any negative spectrum must take the form of discrete eigenvalues. By Theorem~\ref{th: Wu relation}, there are no other eigenvalues and all eigenvalues are simple.

Our variational principle imposes the constraints that each $\lambda_n$ is an eigenvalue of $L_u$.  The structure of eigenfunctions will also play a central role in this section. For any  eigenvalue $\la$, we say that a corresponding eigenfunction $f_\lambda$ is \emph{normalized} if it satisfies
\begin{align}\label{normalization}
\|f_\la\|_{L^2}=1\qquad \text{and} \qquad \langle u, f_\la\rangle>0.
\end{align}
Theorem~\ref{th: Wu relation} guarantees that such a normalized eigenfunction exists and that it is unique.  We will write $f_n$ for the normalized eigenfunction associated to the prescribed eigenvalue $\lambda_n$.

To continue, let  \( P_\text{d}\) and \( P_{\text{e}} \) denote the projections onto the discrete and essential spectrum of \( L_u \), respectively.  For $\kappa\geq C_s(1+\|u\|_{H^s})^{\frac2{1+2s}}$, we may use \eqref{4.1} to deduce that
\begin{align}\label{11:35}
\beta(\kappa;u)
&=\ps{ (L_u+\kappa)^{-1/2} u_+}{ (P_{\text{e}}+P_{\text{d}})(L_u+\kappa)^{-1/2}u_+}\notag\\
&=\ps{(L_u+\kappa)^{-1/2} u_+}{P_{\text{e}}(L_u+\kappa)^{-1/2}u_+}\,+ \sum_{\la\in \sigma_{\text{d}}(L_u)}\frac{|\langle u, f_\la\rangle|^2}{\la+\kappa}\notag\\
&\geq \sum_{n=1}^N \frac{|\langle u, f_n \rangle|^2}{\lambda_n + \kappa} 
= \sum_{n=1}^N \frac{2\pi |\lambda_n|}{\lambda_n + \kappa},
\end{align}
where the last equality follows from Theorem~\ref{th: Wu relation}.  This proves the claim \eqref{1:45} in Theorem~\ref{Th: variational problem}.

Next, we observe that all multisoliton profiles $u=Q_{\Lambda, \vec c}$ with $\vec c \in \R^N$ yield equality in \eqref{11:35} (and so in \eqref{1:45}).  Indeed, if $u=Q_{\Lambda, \vec c}$ for some $\vec c\in \R^N$, Theorem~\ref{T: Sun} guarantees that
$$
\sigma_\text{d}(L_u) = \{\la_1, \ldots, \la_N\} \qtq{and} u_+\in \spn\{f_1, \ldots, f_N\}.
$$
In particular, $P_{\text{e}}u_+=0$ and the computation above simplifies to
\begin{align}\label{u is optimal}
\beta(\kappa;u)=\ps{u_+}{P_{\text{d}}(L_u+\kappa)^{-1}u_+}= \sum_{n=1}^N \frac{2\pi |\lambda_n|}{\lambda_n + \kappa}.
\end{align}

It remains to prove that \( N \)-multisolitons are the \emph{only} optimizers.  Our approach will be to show that any optimizer satisfies the three conditions of Theorem~\ref{T: Sun}.

In order for equality to hold in \eqref{11:35}, there can be no contribution from the essential spectrum.  This forces $P_\text{e} (L_u+\kappa)^{-1/2} u_+ =0$.   Consequently, the function $(L_u+\kappa)^{-1/2} u_+ \in L^2(\R)$ belongs to the subspace spanned by the eigenfunctions of $L_u$ corresponding to negative eigenvalues. Moreover, the Wu relation \eqref{Wu relation Hs} shows that every such eigenvalue contributes and consequently, for equality to hold in \eqref{11:35}, there can be no eigenvalues other than the prescribed $\lambda_n$.  This verifies (ii) from Theorem~\ref{T: Sun} and also shows that $(L_u+\kappa)^{-1/2} u_+\in \spn\{f_1, \ldots, f_N\}$.  Thus
\begin{align}\label{1:17}
u_+\in \spn\{f_1, \ldots, f_N\}, 
\end{align}
which is condition (iii) from Theorem~\ref{T: Sun}.

It remains to verify that optimizers of the variational problem laid out in Theorem~\ref{Th: variational problem} also satisfy hypothesis (i) from Theorem~\ref{T: Sun}, namely, $\langle x\rangle u\in L^2(\R)$.  The relation \eqref{1:17} compels us to prove strong spatial decay of 
eigenfunctions.  This is the topic of the next proposition, which employs a Combes--Thomas argument.  For traditional differential operators, such arguments yield exponential decay of eigenfunctions.  This is simply not true of the operators $L_u$.  For example, the single soliton profile 
\begin{align*}
Q_{-\frac12, 0} (x) = \tfrac2{1+x^2} \qtq{has normalized eigenfunction} f(x)=\tfrac{1}{\sqrt{\pi}} \tfrac{i}{x+i} .
\end{align*}

This example is instructive for two further reasons.  First, it shows that the restriction $\theta<1-\tfrac1p$ in the next proposition cannot be relaxed.  Second, it shows that we cannot hope to show $\langle x\rangle u\in L^2(\R)$ based solely on the decay of eigenfunctions and \eqref{1:17}; we will need to find further cancellation.  Nevertheless, sharp decay of eigenfunctions is a useful start.  The assumption that $u\in H^1(\R)$ in Proposition~\ref{P:CT} will not ultimately be an impediment; see Proposition~\ref{P:1 done}.

\begin{prop}[Eigenfunction decay] \label{P:CT}
Suppose $u\in H^1(\R)$ and $f$ is an eigenfunction of $L_u$.  Then 
\begin{align*}
(1-ix)^\theta f \in L^p_+(\R) \quad\text{whenever} \quad 1<p<\infty \quad\text{and} \quad 0\leq \theta<1-\tfrac1p.
\end{align*}
\end{prop}

\begin{proof}
Theorem~\ref{th: Wu relation} guarantees that the eigenvalue $\lambda$ is negative.  Correspondingly, given $\eta>0$, we may decompose
\begin{align}\label{decomp}
u = u_1 + u_2 \quad \text{with} \quad  u_1\in C_c(\R) \quad \text{and} \quad \|u_2\|_{L^\infty} < \eta|\lambda|.
\end{align}
If $\eta$ is sufficiently small, we may apply Lemma~\ref{L:inverse} to invert the operator $-i\p_x - \la - \CS u_2$ on $L^p_+(\R)$ via Neumann series.   In this way, the eigenfunction equation $L_u f = \lambda f$ may be rewritten as
\begin{align}\label{944}
f=(-i\p_x - \la - \CS u_2)^{-1} \CS (u_1f).
\end{align}

To move forward, we introduce weights into \eqref{944} to obtain
\begin{align} \label{8:21}
(1-i\rho x)^\theta f=\mathcal{T} (1-i\rho x)^\theta \CS (u_1f),
\end{align}
where $\rho>0$ will be chosen shortly and
$$
\mathcal{T}:=\bigl[(1-i\rho x)^\theta (-i\p_x - \la - \CS u_2)(1-i\rho x)^{-\theta}\bigr]^{-1}.
$$

By applying \eqref{H Lp}, we observe that
\begin{align}\label{8:22}
(1-i\rho x)^\theta \CS (u_1f)\in L^p_+(\R) \quad \text{for any $1<p<\infty$ and $0\leq\theta<1-\tfrac1p$}.
\end{align}
This also relies on $f\in H^{1/2}(\R)\hookrightarrow L^q(\R)$ for $2\leq q<\infty$ and the fact that $u_1$ is bounded and of compact support.  

The role of \eqref{8:22} is to control the right-most term in \eqref{8:21}; this leaves us to  prove that the operator $\mathcal T$ is bounded from $L^p_+(\R)$ to $L^p_+(\R)$. To this end, we expand
\begin{equation}\label{8:23}
\begin{aligned}
(1-i\rho x)^\theta &(-i\p_x - \la - \CS u_2)(1-i\rho x)^{-\theta} \\
&=-i\p_x-\la+ \tfrac{\rho\theta}{1-i\rho x} -(1-i\rho x)^\theta\, \CS (1-i\rho x)^{-\theta}u_2
\end{aligned}
\end{equation}
and then employ \eqref{H Lp} to deduce that  
\begin{equation}\label{8:233}
\begin{aligned}
\bigl\|\tfrac{\rho\theta}{1-i\rho x} -(1-i\rho x)^\theta\, &\CS (1-i\rho x)^{-\theta}u_2\bigr\|_{L^p\to L^p} \\
& \lesssim \bigl\|\tfrac{\rho\theta}{1-i\rho x}\bigr\|_{L^\infty}+ \|u_2\|_{L^\infty} \lesssim \rho + \eta|\lambda| .
\end{aligned}
\end{equation} 

Using Lemma~\ref{L:inverse} and choosing $\rho$ and $\eta$ sufficiently small depending on $\la$ and the implicit constant in \eqref{8:233}, we may ensure that \begin{align}\label{8:24}
\bigl\| (-i\p_x-\la)^{-1} \bigr\|_{L^p_+\to L^p_+}\bigl\|\tfrac{\rho\theta}{1-i\rho x} &-(1-i\rho x)^\theta\, \CS (1-i\rho x)^{-\theta}u_2\bigr\|_{L^p\to L^p}\leq \tfrac12.
\end{align}
Using \eqref{8:23} and a Neumann series expansion, we may write
$$
\mathcal T = \sum_{\ell\geq 0}(-i\p_x-\la)^{-1} \Bigl\{ \Bigl[ \tfrac{\rho\theta}{1-i\rho x} -(1-i\rho x)^\theta\, \CS (1-i\rho x)^{-\theta}u_2\Bigr](-i\p_x-\la)^{-1}\Bigr\}^\ell.
$$
The convergence of this series to a bounded operator from $L^p_+(\R)$ to $L^p_+(\R)$ follows from \eqref{8:24} and Lemma~\ref{L:inverse}. Together with \eqref{8:21} and \eqref{8:22}, this shows that $(1-i\rho x)^\theta f\in L^p_+(\R)$ for $1<p<\infty$, $0\leq\theta<1-\tfrac1p$,  and some $\rho>0$.  This implies the claim of the proposition.
\end{proof} 

We now come to the final result of this section, namely, Propostion~\ref{P:1 done}.  It shows that optimizers satisfy property (i) of Theorem~\ref{T: Sun}, thereby completing the proof of Theorem~\ref{Th: variational problem}.  As we have noted earlier, eigenfunction decay alone does not provide sufficient decay for $u(x)$; an additional cancellation is required.  As we will now see, just such a cancellation can be fabricated from the relations
\begin{equation}\label{key cancel}
u=u_+ + \overline{u_+} = \sum_{n=1}^N\ps{f_n}{u}\big(f_n+\overline{f_n}\big)  \qtq{and} \int {u}\big(f_n-\overline{f_n}\big) \,dx =0,
\end{equation}
which follow from  \eqref{1:17} and the normalization of eigenfunctions \eqref{normalization}.

\begin{prop}\label{P:1 done} 
If $u\in H^s(\R)$ is an optimizer of the variational problem laid out in Theorem~\ref{Th: variational problem}, then $u\in H^1(\R)$ and $\langle x \rangle u\in L^2(\R)$.
\end{prop}

\begin{proof}
Let us begin by proving that $u\in H^1(\R)$.  As $u$ is an optimizer, we already know \eqref{key cancel}.  Moreover, each eigenfunction $f_n$ belongs to $\Hohp(\R)$, the form domain of the operator $L_u$.  Thus, we may conclude that $u\in H^{1/2}(\R)$. 

Using the Sobolev embedding \(H^{1/2}(\R) \hookrightarrow L^4(\R)\) and the boundedness of the Cauchy--Szeg\H{o} projection on $L^2(\R)$, we deduce that $\Pi(u f_n)\in L^2_+(\R)$.  Combining this with the eigenvalue equation 
\begin{align}\label{evn}
-i \partial_x f_n - \lambda_n f_n = \Pi(u f_n),
\end{align} 
we conclude that every $f_n \in H^1_+(\R)$.  In view of \eqref{key cancel}, it follows that  \( u\in H^1(\R) \). 

Now that we know \( u\in H^1(\R) \), we may apply  Proposition~\ref{P:CT} to obtain
\begin{align}\label{2/3 in 4}
(1-ix)^{\frac23} f_n(x) \in L^4_+(\R)   \qtq{for each $n$,\quad and then}
\langle x\rangle ^{\frac23} u(x) \in L^4(\R).
\end{align} 
by applying \eqref{key cancel} once more.

We now come to the crux of the matter, namely, showing that
\begin{align}\label{good sum}
\langle x\rangle (f_n + \overline{f_n})\in L^2(\R)  \qtq{for all} 1\leq n\leq N.
\end{align}
The first relation in \eqref{key cancel} shows that this suffices to complete the proof of the proposition; the second relation in \eqref{key cancel} will play a key role in proving it.

The remainder of the proof is devoted to verifying \eqref{good sum}.  The index $n$ may be fixed in all that follows and so we will drop the subscript $n$ on both the eigenvalue $\lambda$ and the eigenfunction $f$.

Combining the eigenvalue equation  \( -i \partial_x f - \CS (uf) = \lambda f \) and its complex conjugate, we find
\begin{equation}\label{5.18}
\bigl(\lvert -i\partial_x \rvert + \lvert \lambda \rvert \bigr)(f + \overline{f}) = \tfrac{1}{2} u (f+\overline{f})+\tfrac{i}{2}H\bigl[u(f-\overline{f})\bigr].
\end{equation}
Our goal is to now apply \eqref{E:inverse'}; this relies on demonstrating adequate decay of the right-hand side.  The first term is easy; indeed, \eqref{2/3 in 4} shows
\begin{align}\label{uf+bar}
\langle x\rangle^{\frac43}  u (f \pm \overline{f})  \in L^2(\R) .
\end{align}

The slow decay of the kernel associated to the Hilbert transform is an impediment to treating the second term in the same way.  Here the second relation in \eqref{key cancel} is key: the function $g=u(f-\overline{f})$ satisfies $\langle x\rangle^{\frac43} g \in L^2(\R)$ and $\int g =0$; correspondingly,
\begin{align*}
[ x,  H ]  g = \pv \ \frac1{\pi} \int_\R\frac{x-y}{x-y}\, g(y) \, dy = 0.
\end{align*}
From this, \eqref{uf+bar}, and the case $p=2$, $\theta=\frac13$ of \eqref{H Lp}, we deduce
\begin{equation}\label{E:Huf}\begin{aligned}
\bigl\|  (1-ix)^{\frac43} H u(f-\overline{f}) \bigr\|_{L^2} &= \bigl\|  (1-ix)^{\frac13}  H (1-ix) u(f-\overline{f}) \bigr\|_{L^2} \\
&\lesssim \bigl\|  (1-ix)^{\frac43} u(f-\overline{f}) \bigr\|_{L^2} < \infty.
\end{aligned}\end{equation}

Combining \eqref{uf+bar} and \eqref{E:Huf} to bound RHS\eqref{5.18} and applying \eqref{E:inverse'} we obtain \eqref{good sum} and so complete the proof of the proposition.
\end{proof}


\section{Concentration Compactness Principle for (BO)}\label{S:CC}
The goal of this section is to develop a concentration compactness principle adapted to our variational characterization of multisolitons. 

\begin{theorem}[Concentration Compactness Principle]\label{T:CC}
Fix $-\frac12<\sigma<s<0$, and let \(\{u_n\}_{n\geq 1}\) be a bounded sequence in $H^s(\R)$. Passing to a subsequence, there exist \(J^* \in \{0\} \cup\mathbb{N} \cup \{+\infty\}\), non-zero profiles  \(U^j \in H^{s}(\R)\), and positions \(x_n^j \in  \mathbb{R}\), such that for each finite \( 0 \leq J \leq J^* \) we have the decomposition
\begin{equation}\label{profile decomp}
u_n(x) = \sum_{j=1}^{J} U^j(x - x_n^j) + r_n^J
\end{equation}
with the following properties: for each fixed $\kappa\geq C_s( 1+\sup_n2\|u_n\|_{H^s})^{\frac2{1+2s}}$,
\begin{gather}
r_n^J(\,\cdot\,+x_n^j)\rightharpoonup 0\quad \text{ in $H^s(\R)$ for each $1\leq j\leq J$}, \label{rnj w-limit 0}    \\ 
\lim_{J \to J^*} \lim_{n\to \infty} \sup_{y \in \mathbb{R}} \|\langle x - y \rangle^{-2} r_n^J \|_{H^{\sigma}} = 0, \label{eq:remainder_decay_weak_norm}  \\
\lim_{J\to J^*}\lim_{n\to \infty}\, \|u_n\|_{H^s}^2-\|r_n^J\|_{H^s}^2-\sum_{j=1}^J\|U^j\|_{H^s}^2=0, \label{eq: Hs decoup}\\
\lim_{J\to J^*} \lim_{n\to\infty} \, \beta(\kappa;u_n)-\sum_{j=1}^{J} \beta(\kappa;U^j)-\beta(\kappa;r_n^J)= 0, \label{eq: beta decoup}\\
\lim_{n \to \infty} |x_n^j - x_n^\ell| = \infty \quad \text{for } j \neq \ell.  \label{xnj orthogonality}
\end{gather}
Moreover,
\begin{equation} \label{rnJ quadratic}  
\lim_{J\to J^*} \lim_{n\to \infty} \, \Bigl| \beta(\kappa;r_n^J) - \int_0^{\infty} \frac{|\widehat{r_n^J}(\xi)|^2}{\xi+\kappa} \,d\xi \Bigr| = 0 .     
\end{equation}
\end{theorem}

We prove Theorem~\ref{T:CC} by induction employing the following inverse inequality:

\begin{prop} \label{Prop: extraction_profil}
Fix $-\frac12<\sigma<s<0$ and suppose \( u_n \in H^s(\R) \) satisfy
\begin{align}\label{Meps}
\|u_n\|_{H^{s}} \leq M \quad \text{and} \quad \sup_{y \in \mathbb{R}} \, \bigl\|\langle x - y \rangle^{-2} u_n \bigr\|_{H^{\sigma}} > \varepsilon
\end{align}
for all $n\geq 1$.  Passing to a subsequence, there exist \( x_n\in \mathbb{R} \) and a non-zero function \( U \in H^s(\R) \) so that
\begin{align}\label{vn to U}
u_n(\cdot + x _n) \rightharpoonup U \quad \text{in } H^s(\R)
\end{align}
and
\begin{align}\label{U lb}
\| U \|_{H^s} \gtrsim \eps \bigl( \tfrac{\eps}{M} \bigr)^{\frac{1 - 2s}{2(s - \sigma)}} .
\end{align}
Moreover, \eqref{vn to U} guarantees the following asymptotic decoupling properties:
\begin{align}
&\lim_{n\to\infty} \|u_n\|_{H^s}^2 - \|U\|_{H^s}^2 - \|u_n(\,\cdot\,+x_n)-U\|_{H^s}^2 = 0, \label{decoup Hs}\\
&\lim_{n\to\infty} \beta(\kappa;u_n)-\beta(\kappa;U)-\beta(\kappa;u_n(\cdot+x_n)-U)= 0, \label{decoup beta}
\end{align}
for $\kappa\geq C_s( 1+2M)^{\frac2{1+2s}}$.
\end{prop}

\begin{proof}
Decomposing into high and low frequencies with respect to a cutoff $N\in 2^\N$ to be determined later, we first use Lemma~\ref{L:multiplier} and Bernstein inequalities to bound the contribution of the high frequencies as follows:
\begin{align*}
\sup_{y\in\R}\|\langle x-y\rangle^{-2}P_{> N} u_{n}\|_{H^{\sigma}}
&\lesssim \|P_{> N} u_{n}\|_{H^{\sigma}}\lesssim N^{\sigma-s}\|P_{> N}u_{n}\|_{H^{s}}\lesssim N^{\sigma-s} M. 
    \end{align*}
Thus, choosing 
$$
N =C \bigl( \tfrac M\eps\bigr)^{\frac{1}{s-\sigma}}
$$
for a constant $C$ that is large enough, we get
$$
\sup_{y\in\R}\|\langle x-y\rangle^{-2}P_{> N} u_{n}\|_{H^{\sigma}}<\tfrac\eps2.
$$
Hypothesis \eqref{Meps} then yields
\begin{align*}
\tfrac\eps2<\sup_{y\in\R}\| \langle x-y\rangle^{-2}P_{\leq N}u_{n}\|_{H^{\sigma}}\leq \sup_{y\in\R}\|\langle x-y\rangle^{-2}P_{\leq N}u_{n}\|_{L^2} \lesssim \|P_{\leq N}u_{n}\|_{L^\infty}.
\end{align*}
Thus, there exist positions $x_n\in \R$ so that 
\begin{align}\label{xn}
\bigl|P_{\leq N}u_{n}\bigr|(x_n)\gtrsim \eps.
\end{align}

As $u_n$ is bounded in $H^s(\R)$, the Banach--Alaoglu Theorem applies.  Thus, we may pass to a subsequence to ensure that \eqref{vn to U} holds, that is, 
\begin{align}\label{weak convg}
v_n:= u_n(\,\cdot+x_n)\rightharpoonup U \quad\text{in} \quad H^{s}(\R).
\end{align}

The linear functional $u\mapsto [P_{\leq N} u](0)$ is $H^s$-continuous; indeed,
$$
|P_{\leq N}u(0)| \leq N^{\frac12-s} \| u \|_{H^s}.
$$
In this way, \eqref{xn} and \eqref{weak convg} imply
\begin{equation}
\| U \|_{H^s} \gtrsim N^{s-\frac12} |P_{\leq N}U(0)|\gtrsim N^{s-\frac12} \eps. 
\end{equation}
Substituting our choice of $N$, we obtain \eqref{U lb}.

We now turn to the asymptotic decoupling statements. Recalling \eqref{weak convg}, we see that $v_n-U \rightharpoonup 0$ in $H^s(\R)$ as $n\to\infty$ and consequently,
\begin{align*}
\|u_n\|_{H^s}^2 - \|U\|_{H^s}^2 - \|u_n(\,\cdot\,+x_n)-U\|_{H^s}^2 = 2 \Re \, \bigl\langle U, v_n- U\bigr\rangle_{H^s} \to 0
\end{align*}
as $n\to\infty$.  This proves \eqref{decoup Hs}.

It remains to prove the decoupling property \eqref{decoup beta}. Passing to a subsequence if necessary, \eqref{Meps} and \eqref{decoup Hs} ensure that 
\begin{align}\label{2:30}
\|U\|_{H^s}\leq 2M \quad \text{and}\quad \sup_n\, \| v_n-U\|_{H^s} \leq 2M.
\end{align}
In particular, $\beta(\kappa;u_n)=\beta(\kappa;v_n)$, $\beta(\kappa;U)$, and $\beta(\kappa;v_n-U)$ are all well-defined for $\kappa\geq C_s( 1+2M)^{\frac2{1+2s}}$.
As
\begin{align*}
\beta(\kappa;v_n)  -\beta(\kappa;U)-\beta(\kappa;v_n-U)
&=\bigl\langle v_n,  (L_{v_n} + \kappa)^{-1}  (v_n)_+ \bigr\rangle-\bigl\langle U,(L_{U} + \kappa)^{-1} U_+ \bigr\rangle\\
&\quad- \bigl\langle  v_n-U,  (L_{v_n-U} + \kappa)^{-1} \big(v_n-U\big)_+ \bigr\rangle,
\end{align*}
adding and subtracting 
$$
\bigl\langle  U,  \,(L_{v_n} + \kappa)^{-1} U_+ \bigr\rangle \quad \text{and}  \quad  \bigl\langle  v_n-U,  \,(L_{v_n} + \kappa)^{-1} \big(v_n-U\big)_+ \bigr\rangle
$$
we may decompose
\begin{align*}
\beta(\kappa;v_n)  -\beta(\kappa;U)-\beta(\kappa;v_n-U)= I_1+I_2+I_3,
\end{align*}
where 
\begin{align*}
I_1&:=\bigl\langle v_n,  (L_{v_n} + \kappa)^{-1}  (v_n)_+ \bigr\rangle \! - \!\bigl\langle  U, (L_{v_n} + \kappa)^{-1}  U_+ \bigr\rangle \! - \! \bigl\langle v_n\! - \! U, (L_{v_n} + \kappa)^{-1} \big(v_n \!-\!U\big)_+ \bigr\rangle\\
&=2\Re\, \bigl\langle  U, (L_{v_n} + \kappa)^{-1}   \big(v_n-U\big)_+ \bigr\rangle,\\
I_2&:=\bigl\langle U,\! \bigl[(L_{v_n} + \kappa)^{-1}\!-\!(L_{U} + \kappa)^{-1}  \bigr] U_+ \bigr\rangle \!= \!\bigl\langle (L_{v_n} + \kappa)^{-1} U_+, (v_n\!-\! U)(L_{U} + \kappa)^{-1}  U_+\bigr\rangle\\
I_3&:=\bigl\langle  v_n-U,\bigl[(L_{v_n} + \kappa)^{-1} -(L_{v_n-U} + \kappa)^{-1} \bigr] \big( v_n-U\big)_+\bigr\rangle\\
&\ =\bigl\langle   (L_{v_n} + \kappa)^{-1}\big(v_n-U\big)_+\,,\, U(L_{v_n-U} + \kappa)^{-1}  \big( v_n-U\big)_+ \bigr\rangle.
\end{align*}
We will prove that all three terms converge to zero as $n\to \infty$, starting with $I_1$.

For $0<\nu<\rho\leq 1$ to be chosen later, we decompose
\begin{align}
I_1&= 2\Re\, \bigl\langle  \bigl[1- \tfrac{1}{1-i\rho x}\bigr]U, (L_{v_n} + \kappa)^{-1}   \bigl(v_n-U\bigr)_+ \bigr\rangle \notag\\
&\quad + 2\Re\, \bigl\langle  \tfrac{1}{1-i\rho x} U, (L_{v_n} + \kappa)^{-1}   \tfrac{1}{1-i\nu x}\bigl(v_n-U\bigr)_+ \bigr\rangle \label{d I1}\\
&\quad + 2\Re\, \bigl\langle  \tfrac{1}{1-i\rho x} U, (L_{v_n} + \kappa)^{-1} \bigl[1- \tfrac{1}{1-i\nu x}\bigr]\bigl(v_n-U\bigr)_+ \bigr\rangle. \notag
\end{align}
Using \eqref{2:30}, the fact that $H^{r+1}(\R)\subset H^{-r}(\R)$ for any $|r|<\frac12$, followed by \eqref{resolvent bound}, we may bound
\begin{align*}
\bigl|\bigl\langle  \bigl[1- \tfrac{1}{1-i\rho x}\bigr]U &, (L_{v_n} + \kappa)^{-1}   \bigl(v_n-U\bigr)_+ \bigr\rangle\bigr|\\
&\lesssim \bigl\|\bigl[1- \tfrac{1}{1-i\rho x}\bigr]U\bigr\|_{H^s} \bigl\|(L_{v_n} + \kappa)^{-1}   \bigl(v_n-U\bigr)_+\bigr\|_{H^{-s}}\\
&\lesssim  \bigl\|\bigl[1- \tfrac{1}{1-i\rho x}\bigr]U\bigr\|_{H^s}  \bigl\|(L_{v_n} + \kappa)^{-1} \bigr\|_{H^s\to H^{s+1}}\|v_n-U\|_{H^{s}}\\
&\lesssim M\bigl\|\bigl[1- \tfrac{1}{1-i\rho x}\bigr]U\bigr\|_{H^s},
\end{align*}
uniformly in $n\geq 1$.  To continue, we use Lemma~\ref{L:multiplier} and the fact that 
$$
\bigl[\tfrac{1}{1-i\rho x} \bigr]\widehat{\phantom{I} } (\xi) = \tfrac{\sqrt{2\pi}}{i\rho}e^{-\xi/\rho}1_{\xi\geq 0}
$$
to estimate
\begin{align*}
\bigl\|\bigl[1- \tfrac{1}{1-i\rho x}\bigr]U\bigr\|_{H^s}
&\leq \bigl\|1- \tfrac{1}{1-i\rho x} \bigr\|_{W^{1,\infty}} \|U-\phi\|_{H^s} + \bigl\|\bigl[1- \tfrac{1}{1-i\rho x}\bigr]\phi\bigr\|_{H^s} \\
&\lesssim \|U-\phi\|_{H^s} + \Bigl\| \langle \xi\rangle^s \int_0^\infty \bigl[\widehat\phi(\xi) -\widehat \phi(\xi-\rho \eta)\bigr] e^{-\eta}\, d\eta\Bigr\|_{L^2_\xi}.
\end{align*}
Approximating $U$ by Schwartz functions $\phi$ and using the fact that for such functions we may employ the Dominated Convergence Theorem to conclude that
$$
 \Bigl\| \langle \xi\rangle^s \int_0^\infty \bigl[\widehat\phi(\xi) -\widehat \phi(\xi-\rho \eta)\bigr] e^{-\eta}\, d\eta\Bigr\|_{L^2_\xi} \to 0 \quad\text{as} \quad \rho\to0,
$$
we deduce that 
\begin{align}\label{localized U}
\bigl\|\bigl[1- \tfrac{1}{1-i\rho x}\bigr]U\bigr\|_{H^s} \to 0 \quad\text{as} \quad \rho\to 0.
\end{align}
Consequently, for any $\delta>0$ there exists $\rho(\delta)>0$ so that 
\begin{align}\label{3:56}
\bigl|\bigl\langle  \bigl[1- \tfrac{1}{1-i\rho x}\bigr]U, (L_{v_n} + \kappa)^{-1}   \bigl(v_n-U\bigr)_+ \bigr\rangle\bigr|\leq \tfrac{\delta}{100} \quad\text{for all} \quad \rho<\rho(\delta),
\end{align}
uniformly for $n\geq 1$.

To control the second term in the decomposition \eqref{d I1} of $I_1$, we use \eqref{weak convg} and Lemma~\ref{L:Riesz} to deduce that
\begin{align}\label{strong convg}
\bigl\| \tfrac{1}{1-i\nu x}\bigl(v_n-U\bigr)_+\bigr\|_{H^\sigma} \to 0 \quad\text{as}\quad n\to \infty.
\end{align}
Thus,  using Lemma~\ref{L:multiplier} and \eqref{resolvent bound}, we may bound
\begin{align}\label{3:57}
\bigl|\bigl\langle  \tfrac{1}{1-i\rho x}& U, \,(L_{v_n} + \kappa)^{-1}   \tfrac{1}{1-i\nu x}\bigl(v_n-U\bigr)_+ \bigr\rangle\bigr|\notag\\
&\lesssim \bigl\|  \tfrac{1}{1-i\rho x} U \bigr\|_{H^\sigma} \bigl\| (L_{v_n} + \kappa)^{-1} \tfrac{1}{1-i\nu x}\bigl(v_n-U\bigr)_+\bigr\|_{H^{-\sigma}}\notag\\[-4mm]
&\lesssim \bigl\|  \tfrac{1}{1-i\rho x}\bigr\|_{W^{1,\infty}} \|U\|_{H^\sigma} \bigl\|(L_{v_n} + \kappa)^{-1} \bigr\|_{H^\sigma\to H^{\sigma+1}}\bigl\| \tfrac{1}{1-i\nu x}\bigl(v_n-U\bigr)_+\bigr\|_{H^\sigma} \leq \tfrac{\delta}{100},
\end{align}
whenever $n\geq n(\delta,\nu)\gg1$.

Arguing similarly and using also Proposition~\ref{P:locality} with $0<\theta<\frac12$ and $0\leq \rho\leq 1$, we estimate the third term in the decomposition of $I_1$ as follows:
\begin{align}\label{3:58}
\bigl|\bigl\langle  \tfrac{1}{1-i\rho x} &\, U, (L_{v_n} + \kappa)^{-1} \bigl[1- \tfrac{1}{1-i\nu x}\bigr]\bigl(v_n-U\bigr)_+ \bigr\rangle\bigr|\notag\\
&\lesssim \bigl\|  \bigl(\tfrac{1}{1-i\rho x}\bigr)^{1-\theta}\bigr\|_{W^{1,\infty}}\|U\|_{H^s} \bigl\|(1-i\rho x)^{-\theta}(L_{v_n} + \kappa)^{-1}(1-i\rho x)^{\theta}\|_{H^s\to H^{s+1}}\notag\\
&\qquad\times \bigl\|  \bigl(\tfrac{1}{1-i\rho x}\bigr)^{\theta} \tfrac{\nu x}{1-i\nu x}\bigr\|_{W^{1,\infty}} \|v_n-U\|_{H^s}\notag\\
&\lesssim \bigl\| \bigl(\tfrac{1}{1-i\rho x}\bigr)^{\theta} \tfrac{\nu x}{1-i\nu x}\bigr\|_{W^{1,\infty}} M^2 \lesssim \bigl(\tfrac{\nu}{\rho} \bigr)^\theta M^2 \leq\tfrac\delta{100},
\end{align}
uniformly for $n\geq 1$, provided $\nu\leq \nu_0(\rho, \delta)$.

Collecting \eqref{3:56}, \eqref{3:57}, and \eqref{3:58}, we conclude that
$$
|I_1| \leq \tfrac{3\delta}{50} \quad \text{whenever $n\geq n(\delta)\gg1$.}
$$

We now turn to $I_2$, which we decompose as follows:
\begin{align*}
I_2&=\bigl\langle (L_{v_n} + \kappa)^{-1} \bigl[1-\tfrac1{1-i\rho x}\bigr]U_+, (v_n-U)(L_{U} + \kappa)^{-1}  U_+\bigr\rangle\\
&\quad + \bigl\langle (L_{v_n} + \kappa)^{-1} \tfrac1{1-i\rho x} U_+, \tfrac1{1-i\nu x}(v_n-U)(L_{U} + \kappa)^{-1}  U_+\bigr\rangle\\
&\quad + \bigl\langle (L_{v_n} + \kappa)^{-1} \tfrac1{1-i\rho x} U_+, \bigl[1-\tfrac1{1-i\nu x}\bigr](v_n-U)(L_{U} + \kappa)^{-1}  U_+\bigr\rangle.
\end{align*}
Arguing as for \eqref{3:56} and employing \eqref{ke1} and \eqref{localized U}, we may bound
\begin{align*}
\bigl|\bigl\langle (L_{v_n} + \kappa)^{-1}& \bigl[1-\tfrac1{1-i\rho x}\bigr]U_+, (v_n-U)(L_{U} + \kappa)^{-1}  U_+\bigr\rangle\bigr|\\
&\lesssim \bigl\|(L_{v_n} + \kappa)^{-1} \bigl[1-\tfrac1{1-i\rho x}\bigr]U_+\bigr\|_{H^{-s}} \bigl\|(v_n-U)(L_{U} + \kappa)^{-1}  U_+\bigr\|_{H^s}\\
&\lesssim  \bigl\|(L_{v_n} + \kappa)^{-1} \bigr\|_{H^s\to H^{s+1}}\bigl\|\bigl[1-\tfrac1{1-i\rho x}\bigr]U_+\bigr\|_{H^s} \|v_n-U\|_{H^s} \\
&\qquad \times \bigl\|(L_{U} + \kappa)^{-1}\bigr\|_{H^s\to H^{s+1}}\| U\|_{H^s}\\
&\lesssim M^2 \bigl\|\bigl[1-\tfrac1{1-i\rho x}\bigr]U_+\bigr\|_{H^s} \leq \tfrac\delta{100} \quad\text{for all} \quad \rho<\rho(\delta),
\end{align*}
uniformly for $n\geq1$.  Arguing as for \eqref{3:57} and using \eqref{strong convg}, we may bound
\begin{align*}
\bigl| \bigl\langle (L_{v_n} + \kappa)^{-1}& \tfrac1{1-i\rho x} U_+, \tfrac1{1-i\nu x}(v_n-U)(L_{U} + \kappa)^{-1}  U_+\bigr\rangle\bigr|\\
&\lesssim\bigl\|(L_{v_n} + \kappa)^{-1} \bigr\|_{H^\sigma\to H^{\sigma+1}}\bigl\|\tfrac1{1-i\rho x} U_+\bigr\|_{H^\sigma} \bigl\|\tfrac1{1-i\nu x}(v_n-U)\bigr\|_{H^\sigma}\\
&\qquad \times\bigl\|(L_{U} + \kappa)^{-1}\bigr\|_{H^\sigma\to H^{\sigma+1}}\| U\|_{H^\sigma}\\
&\lesssim M^2 \bigl\|\tfrac1{1-i\nu x}(v_n-U)\bigr\|_{H^\sigma}\leq \tfrac{\delta}{100} \quad\text{for all} \quad n\geq n(\delta,\nu).
\end{align*}
Finally, arguing as for \eqref{3:58} we may bound
\begin{align*}
\bigl| \bigl\langle (L_{v_n} &+ \kappa)^{-1} \tfrac1{1-i\rho x} U_+, \bigl[1-\tfrac1{1-i\nu x}\bigr](v_n-U)(L_{U} + \kappa)^{-1}  U_+\bigr\rangle\bigr|\\
&\lesssim \bigl\|\bigl(\tfrac1{1-i\rho x}\bigr)^{1-\theta} U_+\bigr\|_{H^s}  \bigl\|(1-i\rho x)^{\theta}(L_{v_n} + \kappa)^{-1}(1-i\rho x)^{-\theta}\|_{H^s\to H^{s+1}}\\
&\qquad\times \bigl\| \bigl(\tfrac{1}{1-i\rho x}\bigr)^{\theta} \tfrac{\nu x}{1-i\nu x}\bigr\|_{W^{1,\infty}} \|v_n-U\|_{H^s}\bigl\|(L_{U} + \kappa)^{-1}\bigr\|_{H^s\to H^{s+1}}\| U\|_{H^s}\\
&\lesssim \bigl(\tfrac{\nu}{\rho} \bigr)^\theta M^2 \leq\tfrac\delta{100},
\end{align*}
uniformly for $n\geq 1$, provided $\nu\leq \nu_0(\rho, \delta)$.  Collecting these estimates we obtain
$$
|I_2| \leq \tfrac{3\delta}{100}\quad \text{whenever $n\geq n(\delta)\gg1$.}
$$

Finally, we turn to $I_3$, which we further decompose as follows:
\begin{align*}
I_3& = \bigl\langle   (L_{v_n} + \kappa)^{-1}\big(v_n-U\big)_+\,,\,  \bigl[1-\tfrac1{1-i\rho x}\bigr]U(L_{v_n-U} + \kappa)^{-1}  \big( v_n-U\big)_+ \bigr\rangle\\
&\quad+\bigl\langle   (L_{v_n} + \kappa)^{-1} \tfrac1{1-i\nu x}\big(v_n-U\big)_+\,,\,  \tfrac1{1-i\rho x}U(L_{v_n-U} + \kappa)^{-1}  \big( v_n-U\big)_+ \bigr\rangle\\
&\quad+\bigl\langle   (L_{v_n} + \kappa)^{-1} \bigl[1-\tfrac1{1-i\nu x}\bigr]\big(v_n-U\big)_+\,,\,  \tfrac1{1-i\rho x}U(L_{v_n-U} + \kappa)^{-1}  \big( v_n-U\big)_+ \bigr\rangle.
\end{align*}
Arguing as above, we may bound
\begin{align*}
\bigl|\bigl\langle (L_{v_n} + \kappa)^{-1}&\big(v_n-U\big)_+\,,\,  \bigl[1-\tfrac1{1-i\rho x}\bigr]U(L_{v_n-U} + \kappa)^{-1}  \big( v_n-U\big)_+ \bigr\rangle \bigr|\\
&\lesssim \bigl\|(L_{v_n} + \kappa)^{-1} \bigr\|_{H^s\to H^{s+1}} \bigl\|(L_{v_n-U} + \kappa)^{-1}\bigr\|_{H^s\to H^{s+1}}\|v_n-U\|_{H^s}^2\\
&\qquad \times \bigl\|\bigl[1-\tfrac1{1-i\rho x}\bigr]U_+\bigr\|_{H^s} \\
&\lesssim M^2 \bigl\|\bigl[1-\tfrac1{1-i\rho x}\bigr]U_+\bigr\|_{H^s} \leq \tfrac\delta{100} \quad\text{for all} \quad \rho<\rho(\delta) 
\quad\text{and} \quad n\geq 1,\\[2mm]
\bigl| \bigl\langle   (L_{v_n} + \kappa)^{-1}& \tfrac1{1-i\nu x}\big(v_n-U\big)_+\,,\,  \tfrac1{1-i\rho x}U(L_{v_n-U} + \kappa)^{-1}  \big( v_n-U\big)_+ \bigr\rangle\bigr|\\
&\lesssim  \bigl\|(L_{v_n} + \kappa)^{-1} \bigr\|_{H^\sigma\to H^{\sigma+1}} \bigl\|\tfrac1{1-i\nu x}(v_n-U)\bigr\|_{H^\sigma}\bigl\|\tfrac1{1-i\rho x} U_+\bigr\|_{H^\sigma} \\
&\qquad \times \bigl\|(L_{v_n-U} + \kappa)^{-1}\bigr\|_{H^\sigma\to H^{\sigma+1}}\|v_n-U\|_{H^\sigma}\\
&\lesssim M^2 \bigl\|\tfrac1{1-i\nu x}(v_n-U)\bigr\|_{H^\sigma}\leq \tfrac{\delta}{100} \quad\text{for all} \quad n\geq n(\delta,\nu),\\[2mm]
\bigl|\bigl\langle   (L_{v_n} + \kappa)^{-1} &\bigl[1-\tfrac1{1-i\nu x}\bigr]\big(v_n-U\big)_+\,,\,  \tfrac1{1-i\rho x}U(L_{v_n-U} + \kappa)^{-1}  \big( v_n-U\big)_+ \bigr\rangle \bigr|\\
&\lesssim \bigl\|(1-i\rho x)^{-\theta}(L_{v_n} + \kappa)^{-1}(1-i\rho x)^{\theta}\|_{H^s\to H^{s+1}}\bigl\| \bigl(\tfrac{1}{1-i\rho x}\bigr)^{\theta} \tfrac{\nu x}{1-i\nu x}\bigr\|_{W^{1,\infty}}\\
&\qquad\times \bigl\|\bigl(\tfrac1{1-i\rho x}\bigr)^{1-\theta} U_+\bigr\|_{H^s} \bigl\|(L_{v_n-U} + \kappa)^{-1}\bigr\|_{H^s\to H^{s+1}}\| v_n-U\|_{H^s}^2\\
&\lesssim \bigl(\tfrac{\nu}{\rho} \bigr)^\theta M^2 \leq\tfrac\delta{100}, \quad\text{for all} \quad \nu<\nu(\rho, \delta) \quad\text{and} \quad n\geq 1.
\end{align*}
These estimates lead to 
$$
|I_3| \leq \tfrac{3\delta}{100}\quad \text{whenever $n\geq n(\delta)\gg1$.}
$$

Collecting our bounds on $I_1$, $I_2$, and $I_3$, we find that
\begin{align*}
\bigl|\beta(\kappa;v_n)  -\beta(\kappa;U)-\beta(\kappa;v_n-U)\bigr| <\delta \quad \text{whenever $n\geq n(\delta)\gg1$,}
\end{align*}
which completes the proof of the decoupling statement \eqref{decoup beta}.
\end{proof}

\begin{proof}[Proof of Theorem~\ref{T:CC}]
Fix $\kappa\geq C_s( 1+\sup_n2\|u_n\|_{H^s})^{\frac{2}{1+2s}}$.  We apply Proposition~\ref{Prop: extraction_profil} iteratively, extracting one profile at a time.  To start, we set $r_n^0:=u_n$.  Now suppose we have a decomposition up to level $J-1\geq 0$ satisfying \eqref{eq: Hs decoup}, \eqref{eq: beta decoup}, and the following special case of \eqref{rnj w-limit 0}:
$$
r_n^{J-1}(\,\cdot\,+x_n^{J-1})\rightharpoonup 0\quad \text{ in $H^s(\R)$}.
$$
Passing to a subsequence if necessary, we set
\begin{align*}
M_{J-1}:=\lim_{n\to\infty} \|r_n^{J-1}\|_{H^s} \quad \text{and} \quad \eps_{J-1}:=\lim_{n\to \infty} \sup_{y\in \R} \|\langle x - y \rangle^{-2} r_n^{J-1} \|_{H^{\sigma}}.
\end{align*}

If $\eps_{J-1}=0$, we stop and set $J^*=J-1$.  If not, we apply Proposition~\ref{Prop: extraction_profil} to $r_n^{J-1}$.  Passing to a further subsequence in $n$, this yields a non-zero profile $U^{J}\in H^s(\R)$ and positions $x_n^{J}\in\R$ such that
\begin{align}\label{weak lim}
U^{J}(x)=\wlim_{n\to\infty}r_n^{J-1}\bigl(x+x_n^{J}\bigr).
\end{align}

To continue, we define
\begin{align}\label{rnJ}
r_n^{J}(x):=r_n^{J-1}(x)-U^{J} \bigl(x-x_n^{J}\bigr).
\end{align}
Note that this definition ensures that \eqref{rnj w-limit 0} holds with $j=J$, that is,
\begin{align}\label{rnJ weak}
r_n^{J}(\,\cdot\,+x_n^{J})\rightharpoonup 0\quad \text{ in $H^s(\R)$}.
\end{align}
By Proposition~\ref{Prop: extraction_profil}, we also have
\begin{gather*}
\lim_{n\to\infty} \bigl\|r_n^{J-1}\bigr\|_{H^s}^2 - \bigl\|U^{J}\bigr\|_{H^s}^2 - \bigl\|r_n^{J}\bigr\|_{H^s}^2 = 0,\\
\lim_{n\to\infty} \, \beta(\kappa;r_n^{J-1})-\beta(\kappa;U^{J})- \beta(\kappa;r_n^{J})= 0, 
\end{gather*}
which combined with the inductive hypothesis gives \eqref{eq: Hs decoup} and \eqref{eq: beta decoup} at the level $J$.  In particular, we have
\begin{align*}
\|U^J\|_{H^s}\leq M_{J-1} \quad \text{and}\quad M_{J}\leq M_{J-1}
\end{align*}
and so $\beta(\kappa;U^J)$ and $\beta(\kappa;r_n^J)$ are well-defined for $\kappa\geq C_s( 1+2M_0)^{\frac{2}{1+2s}}$.

This completes the inductive process. There are two possible outcomes: If $\eps_{J}=0$ for some finite $J$, we declare $J^*=J$.  If $\eps_{J}>0$ for each finite $J$, we set $J^*=\infty$.  We must now verify the remaining claims about this decomposition, namely, \eqref{eq:remainder_decay_weak_norm}, \eqref{xnj orthogonality},  and \eqref{rnJ quadratic}, as well as the $1\leq j<J$ cases of \eqref{rnj w-limit 0}.

If $J^*<\infty$, then \eqref{eq:remainder_decay_weak_norm} is automatic. If $J^*=\infty$, we argue as follows: by \eqref{U lb},
\begin{align*}
\|U^j\|_{H^s} \gtrsim \varepsilon_{j-1} \bigl( \tfrac{\varepsilon_{j-1}}{M_{j-1}} \bigr)^{\frac{1 - 2s}{2(s - \sigma)}}\gtrsim \varepsilon_{j-1} \bigl( \tfrac{\varepsilon_{j-1}}{M_0} \bigr)^{\frac{1 - 2s}{2(s - \sigma)}}
\end{align*}
for all $j\geq 1$.  Together with \eqref{eq: Hs decoup} this shows that $\lim_{J\to \infty}\eps_J= 0$ and so \eqref{eq:remainder_decay_weak_norm} also holds when $J^*=\infty$.

We now turn to the asymptotic orthogonality statement \eqref{xnj orthogonality}.  We proceed by contradiction, choosing \((j, \ell) \) that satisfies the following: $j<\ell$, \eqref{xnj orthogonality} fails for $(j,\ell)$,  and \eqref{xnj orthogonality} holds for all pairs $(j,m)$ with $j<m<\ell$.  Passing to a subsequence, we may assume
\begin{align}\label{cg}
\lim_{n\to \infty}\bigl(x_n^j-x_n^\ell\bigr)= x_0 \in \R.
\end{align}
From the inductive relation
\begin{align*}
r_n^{\ell-1}=r_n^j-\sum_{m=j+1}^{\ell-1}U^m(\cdot-x_n^m),
\end{align*}
we get
\begin{align}\label{tp}
U^\ell(x)&=\wlim_{n\to\infty}r_n^{\ell-1}(x+x_n^\ell)\notag\\
&=\wlim_{n\to\infty}r_n^j(x+x_n^\ell)- \sum_{m=j+1}^{\ell-1} \wlim_{n\to \infty}U^m(x+x_n^\ell-x_n^m),
\end{align}
where the weak limits are in the $H^s$ topology.  From \eqref{cg} and \eqref{rnJ weak}, we see that  the first limit on the right-hand side of \eqref{tp} is zero.
That the remaining limits are zero follows from our assumption that \eqref{xnj orthogonality} holds for all pairs $(j,m)$ with $j<m<\ell$ together with \eqref{cg}.  Thus \eqref{tp} yields $U^\ell=0$, which contradicts the nontriviality of $U^\ell$.  This completes the proof of \eqref{xnj orthogonality}.

Claim \eqref{rnj w-limit 0} for $1\leq j<J$ follows from the decomposition
$$
r_n^{J}(x+x_n^j)=r_n^j(x+x_n^j) -\sum_{m=j+1}^{J}U^m(x+x_n^j-x_n^m),
$$
together with \eqref{rnJ weak} and \eqref{xnj orthogonality}.

It remains to prove the claim \eqref{rnJ quadratic}. To this end, we bound
\begin{align}\label{431}
\Bigl|\beta(\kappa;r_n^J) - \int_0^{\infty} \frac{|\widehat{r_n^J}(\xi)|^2}{\xi+\kappa} \,d\xi \Bigr|
&=\bigl| \bigl\langle r_n^J, \bigl[(L_{r_n^J}+\kappa)^{-1}-(-i\partial+\kappa)^{-1}\bigr](r_n^J)_+\bigr\rangle\bigr|\notag\\
&\lesssim \|r_n^J\|_{H^\sigma}^2\bigl\|(L_{r_n^J}+\kappa)^{-1}-(-i\partial+\kappa)^{-1}\bigr\|_{H^\sigma\to H^{\sigma+1}}\notag\\
&\lesssim M_0^2 \bigl\|(L_{r_n^J}+\kappa)^{-1}\CS \, r_n^J(-i\partial+\kappa)^{-1}\bigr\|_{H^\sigma\to H^{\sigma+1}}\notag\\
&\lesssim M_0^2 \bigl\|(L_{r_n^J}+\kappa)^{-1}\bigr\|_{H^\sigma\to H^{\sigma+1}}\bigl\|r_n^J(-i\partial+\kappa)^{-1}\bigr\|_{H^\sigma\to H^{\sigma}}\notag\\
&\lesssim M_0^2 \bigl\|r_n^J(-i\partial+\kappa)^{-1}\bigr\|_{H^\sigma\to H^{\sigma}}.
\end{align}

For $f\in H^\sigma(\R)$ we use Lemma~\ref{L:multiplier} and \eqref{ke1} to obtain
\begin{align*}
\bigl\|r_n^J(-i\partial+\kappa)^{-1}f\bigr\|_{H^{\sigma}}^2
&\simeq\int_\R \bigl\|\langle x-y\rangle^{-4}r_n^J(-i\partial+\kappa)^{-1}f\bigr\|^2_{{H^{\sigma}}}\,dy \\
&\lesssim \int_\R \bigl\|\langle x-y\rangle^{-2}r_n^J\bigr\|_{H^{\sigma}}^2 \bigl\|\langle x-y\rangle^{-2}(-i\partial+\kappa)^{-1}f\bigr\|^2_{H^{\sigma+1}}\,dy \\
&\lesssim \sup_{y\in\R}\,\bigl\|\langle x-y\rangle^{-2}r_n^J\bigr\|_{H^{\sigma}}^2
	\int_\R\bigl\|\langle x-y\rangle^{-2}(-i\partial+\kappa)^{-1}f\bigr\|^2_{H^{{\sigma}+1}}\,dy \\
&\lesssim \sup_{y\in\R}\,\bigl\|\langle x-y\rangle^{-2}r_n^J\bigr\|_{H^{\sigma}}^2
	\bigl\|(-i\partial+\kappa)^{-1}f\bigr\|^2_{H^{{\sigma}+1}} 
\end{align*}
and consequently,
\begin{align}\label{432}
\bigl\|r_n^J(-i\partial+\kappa)^{-1}\bigr\|_{H^\sigma\to H^{\sigma}}  \lesssim \sup_{y\in\R}\,\bigl\|\langle x-y\rangle^{-2}r_n^J\bigr\|_{H^{\sigma}}^2.
\end{align}

Combining \eqref{431} and \eqref{432} we find
$$
\Bigl|\beta(\kappa;r_n^J) - \int_0^{\infty} \frac{|\widehat{r_n^J}(\xi)|^2}{\xi+\kappa} \,d\xi \Bigr|\lesssim \sup_{y\in\R}\,\bigl\|\langle x-y\rangle^{-2}r_n^J\bigr\|_{H^{\sigma}}^2,
$$
which converges to zero as $n\to \infty$ and $J\to J^*$ in view of \eqref{eq:remainder_decay_weak_norm}.  This completes the proof of \eqref{rnJ quadratic} and so the proof of Theorem~\ref{T:CC}.
\end{proof}


\section{Orbital Stability of (BO)}\label{S:OStab}

This section is dedicated to the proof of Theorem~\ref{Th: Orbital Stability} when $-\frac12<s<0$.  The key object will be the generating function
\begin{equation}\label{beta of z}
\beta(z;u) := \langle u_+ ,\, (L_u+z)^{-1} u_+\rangle.
\end{equation}
Although we have only considered spectral parameters $z\geq \kappa_0(u)$ thus far, we will need to allow more general complex $z$ in this section.  Indeed, this extension is the topic of our first result:

\begin{prop}\label{P:beta of z}
Fix $-\frac12<s<0$.  For \( u \in H^s(\mathbb{R}) \), the generating function \eqref{beta of z} extends mero\-morphically to $\C\setminus(-\infty,0]$ and satisfies the bound
\begin{equation}\label{beta bnd}
|\beta(z;u)| \lesssim  \left(1 + \frac{\kappa_0(u) + |z|}{\dist(z,\sigma(-L_u))}  \right) \|u\|_{H^s_{\kappa_0(u)}}^2 .
\end{equation}
Moreover, this function admits the representation
\begin{equation}\label{measure in Hs}
\beta(z;u) = \int \frac{d\mu_u(E)}{E + z}
\end{equation}
with a positive measure \( \mu_u \) satisfying
\begin{align}\label{measure finite}
    \operatorname{supp}(\mu_u)\subseteq\sigma(L_u)\subseteq(-\kappa_0(u),\infty)  \qtq{and}  \int \frac{d\mu_u(E)}{E + \kappa_0(u)}<\infty .
\end{align}
Lastly, if $\lambda\in\sigma_d(L_u)$, then $|\lambda|$ is a simple pole of $z\mapsto \beta(z;u)$ with residue
\begin{equation}\label{measure(mass point)}
    \mu_u(\{\lambda\})=2\pi|\lambda|.  
\end{equation}
\end{prop}

\begin{proof}
We  abbreviate $\kappa=\kappa_0(u)$.  For each $z\notin\sigma(-L_u)$, we may write
\begin{align*}
(L_u+z)^{-1} = (L_u+\kappa)^{-1/2} \bigl[ (L_u+\kappa) ( L_u+z )^{-1} \bigr] (L_u+\kappa)^{-1/2} .
\end{align*}
Defining $w:=(L_u+\kappa)^{-1/2} u_+ \in L^2_+(\R)$, this allows us to write
\begin{equation}\label{beta w}
\beta(z;u) = \bigl\langle w,  (L_u+\kappa) ( L_u+z )^{-1}  w\bigr\rangle = \int \tfrac{E+\kappa}{E+z} \,d\nu(E) = \int \bigl[1+ \tfrac{\kappa-z}{E+z}\bigr] \,d\nu(E),
\end{equation}
where $\nu$ is the spectral measure associated to $L_u$ and the vector $w\in L^2_+(\R)$ by the spectral theorem.  This theorem guarantees that $\nu$ is a finite positive measure supported by $\sigma(L_u)$. It also satisfies
\begin{equation}\label{nu size}
\nu (\R) = \| w \|_{L^2}^2 \lesssim \|u\|_{H^s_\kappa}^2 .
\end{equation}

Any of the representations in \eqref{beta w} allow us to see that $\beta(z;u)$ is meromorphic in $\C\setminus(-\infty,0]$.
Combining the last representation with \eqref{nu size} gives a simple proof of \eqref{beta bnd}.  After defining $d\mu_u(E) = (E+\kappa)d\nu(E)$, the claims \eqref{measure in Hs} and \eqref{measure finite} follow easily from our analysis.

Only \eqref{measure(mass point)} remains.  By Theorem~\ref{th: Wu relation}, we know that eigenvalues are simple.  If $f$ is the normalized eigenfunction associated to the eigenvalue $\lambda$, then
\begin{equation}\label{nu size'}
\mu_u(\{\lambda\}) = (\lambda+\kappa) \nu (\{\lambda\}) = (\lambda+\kappa) | \langle f, w\rangle |^2 = | \langle f, u\rangle |^2 = 2\pi |\lambda|,
\end{equation}
where the last step is an application of the Wu relation \eqref{Wu relation Hs}.
\end{proof}

\begin{prop}\label{P:poles}
Fix $-\frac{1}{2}<s<0$ as well as a set $\Lambda$ of negative parameters \( \lambda_1 < \cdots < \lambda_N \).  Suppose \( \{u_n\}_{n\geq 1} \) is a sequence in \( H^s(\mathbb{R}) \) satisfying
\begin{align}\label{4:06}
M:=\sup_{n\geq 1}\, \| u_n \|_{H^s} <\infty
\end{align}
and
\be\label{beta(un) to beta(soliton)}
\lim_{n\to \infty}\, \beta(\kappa;u_n)=\sum_{\lambda\in \Lambda} \tfrac{2\pi |\lambda|}{\lambda + \kappa} \qtq{for each} \kappa\geq \kappa_1:=C_s( 1+2M)^{\frac2{1+2s}}.
\ee
If there exists a sequence \( \{x_n\}_{n\geq 1} \subseteq \mathbb{R} \) such that
\begin{align}\label{4:07}
u_n(\cdot + x_n) \rightharpoonup U \quad \text{in } H^s(\R) \qtq{with}  U \not\equiv 0,
\end{align}
then there exists a non-empty subset $\Lambda_1$ of \(\Lambda\) such that for each $\kappa\geq \kappa_1$,
\be\label{beta[U] beta[vn]}
\beta(\kappa;U)=\sum_{\la\in\Lambda_1}\tfrac{2\pi |\lambda|}{\lambda + \kappa}\ \quad\text{and} \quad \lim_{n\to \infty} \beta\bigl(\kappa; u_n-U(\cdot-x_n)\bigr)= \sum_{\la\in\Lambda\setminus\Lambda_1}\tfrac{2\pi |\lambda|}{\lambda + \kappa}.
\ee
Consequently, by Theorem~\ref{Th: variational problem}, $U$ is a multisoliton profile with parameters $\Lambda_1$.
\end{prop}

\begin{proof}
We write $v_n:= u_n-U(\cdot-x_n)$. In view of \eqref{4:06} and \eqref{4:07}, we have
\begin{align*}
\|U\|_{H^s}\leq M \qtq{and} \sup_{n\geq 1}\|v_n\|_{H^s}\leq 2M.
\end{align*}
Recalling \eqref{small pert}, we see that $\kappa_1\geq \kappa_0(U)$ and $\kappa_1\geq \kappa_0(v_n)$ for all $n\geq 1$ and so,
\begin{align}\label{4:08}
\sigma(L_U)\subset (-\kappa_1, \infty)\qtq{and} \sigma(L_{v_n})\subset (-\kappa_1, \infty).
\end{align}
In particular, $\beta(\kappa; U)$ and $ \beta(\kappa; v_n)$ are well-defined for all $n\geq 1$ and $\kappa\geq \kappa_1$.

Using the asymptotic decoupling property \eqref{decoup beta} together with \eqref{beta(un) to beta(soliton)}, we find 
\begin{equation}\label{eq: beta[rn]+beta[U] to sum}
\beta(\kappa; U) + \lim_{n\to \infty} \beta(\kappa; v_n) = \sum_{\lambda\in \Lambda} \tfrac{2\pi |\lambda|}{\lambda + \kappa} \qtq{for each} \kappa\geq \kappa_1.
\end{equation}

Next, we extend this convergence from $\kappa\geq \kappa_1$ to a much larger set in the complex plane. Specifically, we will show
\begin{equation}\label{complex pain}
\beta(z; U) + \lim_{n\to \infty} \beta(z; v_n) =\sum_{\lambda\in \Lambda} \tfrac{2\pi |\lambda|}{\lambda + z} \qtq{for all} z\in \Omega:=\C\setminus(-\infty, \kappa_1].
\end{equation}
As we will explain, this is a simple consequence of Montel's Theorem.   By Proposition~\ref{P:beta of z}, the functions $\beta(z; U)$ and $\beta(z; {v_n})$ are holomorphic and locally uniformly bounded on $\Omega$.  Therefore, by Montel's Theorem, any subsequence of $\beta(z;v_n)$ admits a further subsequence that converges uniformly on compact subsets of $\Omega$. 
As \eqref{eq: beta[rn]+beta[U] to sum} uniquely identifies such limits, we deduce that $\beta(z; v_n)$ converges (without passing to a subsequence) uniformly for $z$ in compact subsets of $\Omega$.

Now, employing identity \eqref{measure in Hs} of Proposition~\ref{P:beta of z}, we may rewrite \eqref{complex pain} as 
\begin{equation}\label{convergence with measures}
 \int\frac{d\mu_U(E)}{E+z}+\lim_{n\to \infty}\, \int\frac{d\mu_{v_n}(E)}{E+z}=\sum_{\lambda\in \Lambda} \frac{2\pi |\lambda|}{\lambda + z} \qtq{for all} z\in \Omega,
\end{equation}
where \( \mu_U \) and \( \mu_{v_n} \) denote the measures associated with \( L_U \) and \( L_{v_n} \), respectively.

Moving forward, we will use the relation \eqref{convergence with measures} to locate the poles of $\beta(z;U)$ and $\beta(z;v_n)$.  This has two stages: First we demonstrate pole-free regions (for $n$ large) in each of the $N+1$ gaps between $0 < |\lambda_N| < \cdots < |\lambda_1| < \infty$.  This allows us to construct a system of contours surrounding the pole locations.  These contours are used in the second stage to prove that each $\lambda_m$ is either a pole of $\beta(z;U)$ or the limit of poles of $\beta(z;v_n)$.

Due to the sign flip between the poles of $\beta(z;u)$ and $\sigma(L_u)$, it is clearer to work first with the measures.  Given $a\in(-\infty,0)$ with
\begin{equation}\label{good a}
\dist(a, \Lambda \cup \{0\}) > 0, \qtq{set}  \delta =\tfrac12 \min\Bigl\{\dist(a, \Lambda \cup \{0\}), \, \tfrac{\dist(a, \, \Lambda \cup \{0\})^2 }{ 1 + \sum |\lambda_m| }\Bigr\}.
\end{equation}
We will show that for $n$ sufficiently large depending only on $a$, we have
\begin{equation}\label{not in gaps}
\mu_U\bigl([a - \delta, a + \delta]\bigr) = 0  \qtq{and} \mu_{v_n}\bigl([a - \delta, a + \delta]\bigr) = 0. 
\end{equation}

Taking the imaginary part of \eqref{convergence with measures} with $z=-a+i\delta$ and multiplying by $\delta$ gives
\[
\int\frac{\delta^2\,d\mu_U(E)}{(E-a)^2+\delta^2}+\int\frac{\delta^2\,d\mu_{v_n}(E)}{(E-a)^2+\delta^2}\longrightarrow \sum_{\lambda\in \Lambda} \frac{2\pi |\lambda| \delta^2}{(\lambda-a)^2+\delta^2}\quad \text{as } n\to \infty.
\]
For the specific choice of $\delta$ in \eqref{good a}, this yields
\begin{align}
\tfrac12 \mu_U\bigl([a - \delta, a + \delta]\bigr)&+ \tfrac12 \limsup_{n\to \infty}\mu_{v_n}\bigl([a - \delta, a + \delta]\bigr)\notag\\
&\qquad \leq \sum_{\lambda\in \Lambda}\frac{2\pi \delta^2|\la|}{\dist\big(a,\Lambda\cup\{0\}\big)^2+\delta^2}< \pi \delta.\label{4:11}
\end{align}

On the other hand, Proposition~\ref{P:beta of z} shows that 
\begin{equation}\label{mu[mass point]}
\text{either }\,\, \mu_U\bigl([a - \delta, a + \delta]\bigr) = 0 \quad\text{or }\,\,  \mu_U\bigl([a - \delta, a + \delta]\bigr) \geq 2\pi |a + \delta| \geq 2\pi\delta,
\end{equation}
where the last inequality follows from our choice of $\delta$.  The analogous statement for $\mu_{v_n}([a - \delta, a + \delta])$ is also satisfied for each $n\geq 1$.

Combining \eqref{4:11} and \eqref{mu[mass point]}, we deduce that \eqref{not in gaps} holds whenever $n$ is sufficiently large.  By Proposition~\ref{P:beta of z}, this precludes $\beta(z;U)$ from having poles in $[|a|-\delta, |a|+\delta]$, and likewise for $\beta(z;v_n)$ when $n$ is large.

Next we use these pole-free regions to enlarge the domain in which \eqref{complex pain} holds and then to build the system of contours illustrated in Figure~\ref{F:1}.  First, we fix \( a_1, \ldots, a_{N+1} < 0 \) interlacing with the $\lambda_m$, or  equivalently,
\[
0< |a_{N+1}| < |\lambda_N| < |a_{N}| < |\lambda_{N-1}| < \cdots < |a_2| < |\lambda_1| < |a_{1}| < \infty.
\]

As we have just shown, there exist $\delta>0$ and $N_0\in\N$ such that the intervals $[|a_m|-\delta,|a_m|+\delta]$ contain no poles of $\beta(z;U)$ nor of $\beta(z;v_n)$ when $n\geq N_0$.  Equivalently, these intervals contain no points from $\sigma(-L_U)$ or from $\sigma(-L_{v_n})$.  In this way, Proposition~\ref{P:beta of z} shows that  $\beta(z;U)$ and $\beta(z;v_n)$ are both holomorphic and locally uniformly bounded on the larger region 
\[
\Omega_1 := \Omega \cup \bigl(|a_{N+1}|-\delta,|a_{N+1}|+\delta\bigr)\cup \cdots \cup \bigl(|a_{1}|-\delta,|a_{1}|+\delta\bigr).
\]
Repeating the earlier argument based on Montel's Theorem, we deduce that the convergence in \eqref{complex pain} holds uniformly on compact subsets of $\Omega_1$.

As illustrated in Figure~\ref{F:1}, for each $1\leq m \leq N$ we draw a contour \( \Gamma_m \), which passes through the real axis at $|a_{m+1}|$ and $|a_m|$ and encloses a single \( |\lambda_m|\).    As the contours lie in $\Omega_1$, they avoid the poles of \( \beta(z; U) \) and \( \beta(z; v_n) \), at least for $n\geq N_0$.  All contours are oriented anticlockwise.

\begin{figure}[h!]
\centering
\begin{tikzpicture}[scale=1.3]
    \draw[->] (-0.5,0) -- (8.5,0) node[right] {$\operatorname{Re}(z)$};
    \draw[->] (0,-1.3) -- (0,1.3) node[above] {$\operatorname{Im}(z)$};

    \foreach \x/\label in {
        2.0/{|a_N|},
        4.0/{|a_{N-1}|},
        5.5/{|a_2|},
        7.0/{|a_{1}|}
    } {
        \filldraw (\x,0) circle (1pt);
        \node[below=7pt, right=0.1pt] at (\x,0) {\scriptsize $\label$};
    }
    \filldraw (1.0,0) circle (1pt);
    \node[below=7pt, left=0.1pt] at (1.0,0) {\scriptsize $|a_{N+1}|$};
    
    \foreach \x/\label in {
        1.7/{|\lambda_N|},
        3.2/{|\lambda_{N-1}|},
        6.3/{|\lambda_1|}
    } {
        \node[above=0pt] at (\x,0) {\scriptsize $\label$};
        \draw (\x,0) node {$\times$};
    }
    \foreach \x in {
        4.5, 4.75, 5.0
    } {
        \filldraw[blue] (\x,0.25) circle (1pt);
    }
    \draw[blue, thick] (1.5,0) ellipse [x radius=0.5, y radius=0.6]; 
    \draw[blue, thick] (3,0) ellipse [x radius=1.0, y radius=0.7]; 
    \draw[blue, thick] (6.25,0) circle (0.75); 
    \draw[->,blue,thick] (1.53,-0.6) -- (1.531,-0.6);
    \draw[->,blue,thick] (3.00,-0.7) -- (3.01,-0.7);
    \draw[->,blue,thick] (6.29,-0.75) -- (6.2901,-0.75);
    \node[blue] at (1.7,0.75) {\scriptsize $\Gamma_N$};
    \node[blue] at (3.,0.85) {\scriptsize $\Gamma_{N-1}$};
    \node[blue] at (6.25,0.9) {\scriptsize $\Gamma_1$};
\end{tikzpicture}
\caption{For each \( m = 1, \dots, N \), the contour \( \Gamma_m\) is taken to intersect the real axis at \( |a_{m+1}| \) and \( |a_{m}| \), and is chosen so as to enclose  a single pole \( |\lambda_m| \). }
\label{F:1}
\end{figure}
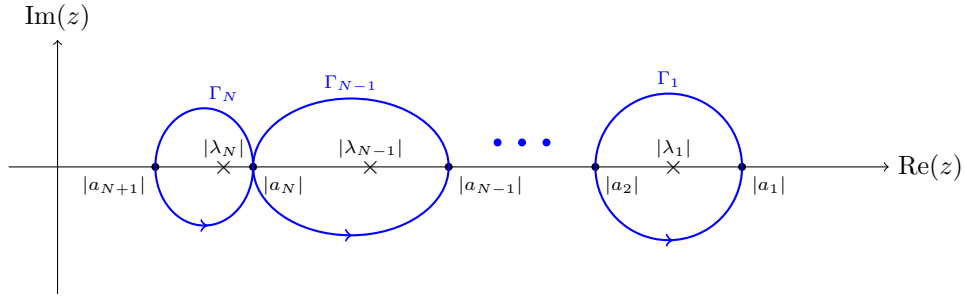

As a proxy for $U$ or $v_n$, let $u\in H^s(\R)$ and let $d\mu_u$ denote the measure associated to it by Proposition~\ref{P:beta of z}.  Assume also that $a_m,a_{m+1}\notin\supp(\mu_u)$.  Given $\ell\in\{-1,0\}$, we apply Fubini, the Cauchy integral formula, and \eqref{measure(mass point)} to obtain
\begin{align*}
\oint_{\Gamma_m} \!z^{\ell}\!  \int \frac{d\mu_u(E)}{E+z} \frac{dz}{2\pi i} &= \int_{a_m}^{a_{m+1}} (-E)^\ell \,d\mu_u(E)
	= \sum_{\lambda\in \sigma_\text{d}(L_u)\cap(a_m,a_{m+1})} 2\pi |\lambda|^{\ell+1} .
\end{align*}
As a special case of this reasoning (or by direct computation), we also have
\begin{align*}
\oint_{\Gamma_m} \!z^{\ell}\, \sum_{j=1}^N \frac{2\pi |\lambda_j |}{\lambda_j + z} \frac{dz}{2\pi i} &= \begin{cases} 2\pi, &\text{if } \ell=-1\\ 2\pi|\lambda_m|, &\text{if } \ell=0. \end{cases}
\end{align*}

Our goal in presenting these computations is to apply them to both sides in  \eqref{convergence with measures}.  Note that we have shown uniform convergence on compact sets, such as  $\Gamma_m$, and this allows us  to exchange the $z$-integral with the limit $n\to\infty$. Proceeding in this manner with $\ell=-1$, we find that for each $1\leq m\leq N$,
\begin{align}
\# \bigl[ \sigma_\text{d}(L_U)\cap(a_m,a_{m+1}) \bigr] +  \# \bigl[ \sigma_\text{d}(L_{v_n})\cap(a_m,a_{m+1}) \bigr] \to 1 \qtq{as} n\to\infty.
\end{align}
For each $m$, we have two mutually exclusive possibilities:  The first is that  $L_U$ has one eigenvalue in $(a_m,a_{m+1})$ and (for $n$ sufficiently large) $L_{v_n}$ has none.  In the second case, $L_{v_n}$ has exactly one eigenvalue (for $n$ sufficiently large) in $(a_m, a_{m+1})$ and $L_U$ has none.  We split $\Lambda=\Lambda_1 \cup \Lambda_2$ accordingly, which is to say, $\lambda_m\in \Lambda_1$ in the first case and $\lambda_m\in \Lambda_2$ in the second case.

Let us now consider the $\ell=0$ case of the computations above.  When $\lambda_m\in \Lambda_1$, we discover that $\lambda_m\in \sigma(L_U)$.  The case $\lambda_m\in \Lambda_2$ is a little more subtle because the sole eigenvalue of $L_{v_n}$ lying in the interval $(a_m,a_{m+1})$ may depend on $n$.  Writing $\lambda_m^{(n)}$ for this eigenvalue, the $\ell=0$ computation shows that $\lambda_m^{(n)} \to \lambda_m$ as $n\to\infty$. 

By Theorem~\ref{Th: variational problem}, we may now infer that
\begin{align}\label{9:07}
\beta(\kappa;U) \geq \sum_{\lambda\in \Lambda_1} \tfrac{2\pi |\lambda|}{\lambda + \kappa}  \qtq{and}
	\lim_{n\to\infty} \beta(\kappa;v_n) \geq \sum_{\lambda\in \Lambda_2} \tfrac{2\pi |\lambda|}{\lambda + \kappa}
\end{align}
for each $\kappa\geq \kappa_1$.  Combining \eqref{eq: beta[rn]+beta[U] to sum} and \eqref{9:07}, we see that equality must occur in both cases in \eqref{9:07}; this settles \eqref{beta[U] beta[vn]}.
Note that this also shows that $\Lambda_1$ is non-empty: The alternative is that 
$$
\beta(\kappa;U)=0 \qtq{for all} \kappa\geq \kappa_1 .
$$
By \eqref{beta est}, this would force $U\equiv 0$, which contradicts the hypothesis \eqref{4:07}.

Finally, the variational characterization of multisolitons from Theorem~\ref{Th: variational problem} shows that \eqref{beta[U] beta[vn]} implies that $U$ is a multisoliton with parameters $\Lambda_1$.
\end{proof}

We are now ready to complete the proof of uniform orbital stability of multisolitons for regularities $-\frac12<s<0$.  The remaining cases will be taken up in the next section.

\begin{proof}[Proof of Theorem~\ref{Th: Orbital Stability} for $-\frac12<s<0$.]
We argue by contradiction. Fix $N\geq 1$ and parameters $\la_1<\cdots<\la_N<0$.  Set $\Lambda=\{\lambda_1, \ldots, \lambda_N\}$ and assume, towards a contradiction,  that there exist $\eps_0>0$, functions \( u_{n,0}\in H^s(\R) \), and times $t_n\in\R$, with
\begin{align}\label{5:09}
\inf_{\vec{c}\in\R^N}\, \bigl\|u_{n,0}-Q_{\Lambda,\vec{c}}\, \bigr\|_{H^s}\longrightarrow 0 \qtq{as} n \to \infty,
\end{align}
but the solutions $u_n$ to \eqref{BO} with initial data $u_n(0)=u_{n,0}$ satisfy
\be\label{contradiction}
\inf_{\vec{c} \in \mathbb{R}^N} \, \bigl\| u_n(t_n) - Q_{\Lambda, \vec{c}} \bigr\|_{H^s} \geq \varepsilon_0 \qtq{for all} n\geq 1.
\ee

From \eqref{Q decay}, \eqref{5:09}, and Proposition~\ref{P:beta}, we see that
\begin{equation}\label{7:bndd}
M:=\sup_{\vec{c}\in\R^N}\,  \|Q_{\Lambda,\vec{c}}\,\|_{H^s} + \sup_{n} \|u_{n,0}\|_{H^s} + \sup_{n} \|u_{n}(t_n)\|_{H^s} <\infty.
\end{equation}
Moreover, \eqref{Q decay} guarantees that $\{ Q_{\Lambda,\vec{c}} : \vec{c}\in\R^N \}$ is $L^2$-bounded and so $H^s$-equicontinuous. Combining this with \eqref{5:09}, we deduce that $\{u_{n,0}: n\in \N\}$ is also $H^s$-equicontinuous. Applying Proposition~\ref{P:beta}, we conclude that 
\begin{equation}\label{7:equi}
\{ u_{n}(t_n) : n\in\N\} \quad\text{is $H^s$-equicontinuous.}
\end{equation}

By \eqref{eq:resolvent_diff_control}, $\beta$ is continuous on \(H^s(\R)\); it is also conserved by the Benjamin--Ono flow.  Thus, combining Theorem~\ref{Th: variational problem} with \eqref{5:09} we obtain
\be\label{beta[un(tn)]=sum}
\beta(\kappa; u_n(t_n))=\beta(\kappa;u_{n,0})\longrightarrow\sum_{\lambda\in \Lambda}\tfrac{2\pi |\la|}{\la+\kappa}  \qtq{as} n \to \infty,
\ee
for all $\kappa\geq C_s(1+M)^{\frac2{1+2s}}$.

From \eqref{7:bndd} we know that the sequence $u_n(t_n)$ is bounded; thus, we may apply our concentration compactness principle, Theorem~\ref{T:CC}.  Passing to a subsequence, we have the decomposition
\be\label{profil_decomp_un(tn)}
u_n(t_n)=\sum_{j=1}^J U^j(\,\cdot-x_n^j)+r_n^J \qtq{for each finite $0\leq J\leq J^*$}
\ee
satisfying all the conclusions of Theorem~\ref{T:CC}. 

Our goal is to prove that there are only finitely many non-zero profiles $U^j$ in this decomposition, that each such $U^j$ is a multisoliton with parameters  $\Lambda_j\subseteq \Lambda$, that the sets $\Lambda_j$ form a partition of $\Lambda$, and that $r_n^J$ converges to $0$ in $H^s(\R)$.  As we will explain, Proposition~\ref{prop:MD} then guarantees the existence of $\vec{c}_n\in \R^N$ such that
\begin{align}\label{12:19}
\| u_n(t_n) - Q_{\Lambda,\vec{c}_n}\, \|_{H^s} \longrightarrow 0 \qtq{as} n\to\infty,
\end{align}
thereby contradicting \eqref{contradiction} and completing the proof.

We start by ruling out the possibility of having no profiles. Assume, towards a contradiction, that $J^*=0$ and so $u_n(t_n)=r_n$. Combining~\eqref{rnJ quadratic} with \eqref{beta[un(tn)]=sum}, we obtain
 \[
\int_0^\infty \frac{|\widehat{r_n}(\xi)|^2}{\xi+\kappa} \, d\xi\ \underset{n\to\infty}\longrightarrow\ \sum_{m=1}^N\frac{2\pi |\la_m|}{\la_m+\kappa}.
 \]
 This is impossible because
\begin{align*}
\kappa\mapsto \kappa\, \int_0^{\infty}\frac{|\widehat{r_n}(\xi)|^2}{\xi+\kappa} \,d\xi\ \,\text{ is increasing,} \qtq{while} \kappa\mapsto\kappa\,\sum_{m=1}^N\frac{2\pi |\la_m|}{\la_m+\kappa}\,\text{ is decreasing}.
\end{align*}
Therefore, we must have that $J^*\geq 1$.

To prove that $J^*<\infty$, we employ Proposition~\ref{P:poles}.  Note that \eqref{beta[U] beta[vn]} allows us to apply Proposition~\ref{P:poles} iteratively to the functions $r_n^J$ defined in \eqref{rnJ}.  In this way, we deduce that for each finite $1\leq j\leq J^*$ there exist disjoint non-empty subsets $\Lambda_j$ of $\Lambda$ such that
\begin{align}\label{beta Uj}
    \beta(\kappa;U^j)= \sum_{\la\in \Lambda_j}\tfrac{2\pi|\la|}{\la+\kappa}.
\end{align}
As $\Lambda$ has cardinality $N$, this shows that $J^*\leq N$.

To see that the sets $\Lambda_j$ partition $\Lambda$, it remains to show that
\begin{align}\label{partition}
\Lambda=\bigcup_{1\leq j\leq J^*} \Lambda_j.
\end{align}
Failure of this would mean that $\Lambda^*:= \Lambda\setminus \bigl( \mkern1mu \cup_{j} \Lambda_j\bigr)$ is non-empty and 
$$
\lim_{n\to \infty} \beta(\kappa; r_n^{J^*}) =\sum_{\lambda\in \Lambda^*}\tfrac{2\pi |\la|}{\la+\kappa}.
$$
Combining this with \eqref{rnJ quadratic}, the competing monotonicities of the functions
$$
\kappa\mapsto \kappa\, \int_0^{\infty}\frac{|\widehat{r_n^{J^*}}(\xi)|^2}{\xi+\kappa} \,d\xi\  \qquad\text{and}\qquad \kappa\mapsto\kappa\,\sum_{\lambda\in \Lambda^*}\frac{2\pi |\la|}{\la+\kappa}
$$
yields a contradiction.

Combining \eqref{eq: beta decoup}, \eqref{beta[un(tn)]=sum}, \eqref{beta Uj}, and \eqref{partition}, we see that
$$
\lim_{n\to\infty}\beta(\kappa; r_n^{J^*}) =0.
$$
Additionally, $r_n^{J^*}$ is $H^s$-equicontinuous; this follows from \eqref{profil_decomp_un(tn)}, \eqref{7:equi}, and the  equicontinuity of the \emph{finite} set $\{U^j:1\leq j\leq J^*\}$.  Thus, Proposition~\ref{P:beta} shows that  $r_n^{J^*}$ converges to zero in $H^s(\R)$.

By \eqref{beta Uj} and the variational characterization of multisolitons given in Theorem~\ref{Th: variational problem}, there exist $\vec{c}_j\in \R^{N_j}$ with $N_j= |\Lambda_j|\in\{1, \ldots, N\}$ so that $U^j=Q_{\Lambda_j, \vec{c}_j}$. Thus, we may rewrite \eqref{profil_decomp_un(tn)} as 
\be\label{un=decomp+rem}
u_n(t_n)=\sum_{j=1}^{J^*} Q_{\Lambda_j, \vec{c}_j}(\cdot-x_n^j)+r_n^{J^*}  \qtq{with} \lim_{n\to \infty}\|r_n^{J^*}\|_{H^s}=0
\ee
and
$$
\lim_{n \to \infty} |x_n^j - x_n^\ell| = \infty \quad \text{for } j \neq \ell. 
$$
Applying Proposition~\ref{prop:MD}, we deduce that there exist $\vec{c}_n\in \R^N$ so that 
$$
\lim_{n\to \infty}\bigl\|Q_{\Lambda, \vec{c}_n} - \sum_{j=1}^{J^*} Q_{\Lambda_j, \vec{c}_j}(\cdot - x_n^j)\bigr\|_{H^s} = 0.
$$
Together with \eqref{un=decomp+rem} this yields \eqref{12:19}, thereby contradicting \eqref{contradiction} and completing the proof of Theorem~\ref{Th: Orbital Stability}.
\end{proof}

\section{Higher regularity}\label{S:HR}

The goal of this section is to prove Theorem~\ref{Th: Orbital Stability} in the cases $0\leq s \leq \frac12$ that were not treated in the previous section.  In fact, we will show how we may deduce these results from those in the last section by using equicontinity and conservation laws.

We find it convenient to prove equicontinuity at higher regularity using a method introduced in \cite{KLV2025}.  The idea is to apply well-chosen operator monotone functions (such as those in Lemma~\ref{L:mono}) to operator inequalities expressing the notion that $L_u$ is a relatively small perturbation of $L_0$ (such as  \eqref{8.1} below).

Taking $s=0$ in \eqref{small pert} and \eqref{small pert'}, we see that there is a constant $C_0\geq 1$ so that
\begin{equation}\label{8.1}\begin{aligned}
\kappa \geq C_0 \|u\|_{L^2}^2 & \implies \| \CS uf \|_{H^{-1/2}_\kappa} \leq \tfrac12 \| f \|_{H^{1/2}_\kappa} \\
& \implies \tfrac12\bigl( L_0 + \kappa ) \leq (L_u + \kappa) \leq \tfrac32\bigl(L_0 + \kappa\bigr),
\end{aligned}\end{equation}
in the sense of quadratic forms.

\begin{lemma}\label{L:mono}
For each $0\leq s \leq \frac12$ and each $\vk\geq 0$, the function
\begin{align*}
	A \mapsto  A (A+\vk)^{2s-1} 
\end{align*}
is operator increasing on positive definite operators.
\end{lemma}

\begin{proof}
When $s=\frac12$, the claim is trivial.  
When $s=0$, we have $A \mapsto 1 - \tfrac{\vk}{A+\vk}$, whose monotonicity follows from that of the operator inverse.  For $0 < s < \frac12$, we exploit the representation
\begin{align*}
	A(A+\vk)^{2s-1} = \tfrac{\sin(2\pi s)}{\pi} \int_0^\infty \Bigl[ 1 - \tfrac{\vk+t}{A+\vk+t} \Bigr] t^{2s-1}\,dt ,
\end{align*}
which shows that the function is a linear combination of operator monotone functions with positive coefficients.  (This is essentially the proof of the easy half of Loewner's characterization of operator monotone functions; see \cite{MR3969971}.) 
\end{proof}

Given $u\in H^s(\R)$, $\kappa$ as in \eqref{8.1}, and $\vk\geq 0$, we define 
\begin{align}\label{B defn}
	B(u,\vk) := \bigl\langle u_+, (L_{u}+\kappa) (L_{u}+\kappa+\vk)^{2s-1} u_+\bigr\rangle.
\end{align}
Observe that the function of $L_u+\kappa$ appearing in \eqref{B defn} is precisely the one shown to be operator monotone in Lemma~\ref{L:mono}.

The quantity $B$ is conserved for any $C^{ }_tH^s$ solution $u(t)$ to \eqref{BO}. To this end, we first note that if $\kappa$ is chosen to satisfy \eqref{8.1} at the initial time, then conservation of $L^2$ guarantees that \eqref{8.1} remains valid throughout the lifespan of the solution.  For smooth solutions, conservation of $B$ follows easily from \eqref{Lax eq} and \eqref{OTP}; this then extends to general solutions by continuity.

\begin{prop}\label{P:s pos}
Fix $0\leq s \leq \frac12$ and let $u(t)$ denote the solution to \eqref{BO} with initial data $u(0)\in H^s(\R)$.  Then 
\begin{align}\label{8:APB}
	\|u(t)\|_{H^s}\lesssim \bigl(1+ \|u(0)\|_{L^2}\bigr)^{2s}\|u(0)\|_{H^s}.
\end{align}
In the case $0\leq s<\frac12$, we also have propagation of equicontinuity: if the set $U$ is $H^s(\R)$-bounded and equicontinuous, then so too is
$$
U^* := \{ u(t) : t\in\R \text{ and } u\text{ is the solution to \eqref{BO} with initial data $u(0)\in U$}\bigr\}.
$$
\end{prop}

\begin{proof}
Choosing $\kappa = 1 + C_0 \|u(0)\|_{L^2}^2$ and using the conservation of $L^2$, we see that condition \eqref{8.1} is satisfied by for all times $t$. By Lemma~\ref{L:mono}, we deduce
\begin{align*}
	2^{-2s} \bigl\langle u_+(t), (L_0+\kappa) (L_0+\kappa+2\vk )^{2s-1} u_+(t)\bigr\rangle &\leq B(u(t),\vk),  \\
	\bigl(\tfrac32\bigr)^{2s} \bigl\langle u_+(0), (L_0+\kappa) (L_0+\kappa+\tfrac23\vk )^{2s-1} u_+(0)\bigr\rangle &\geq B(u(0),\vk) .
\end{align*}
By the conservation of $B$, we then obtain
\begin{align}\label{8:star}
	\int_0^\infty \frac{(\xi+\kappa)|\widehat u(t,\xi)|^2\,d\xi}{(\xi+\kappa+2\vk)^{1-2s}} \lesssim  \int_0^\infty \frac{(\xi+\kappa)|\widehat u(0,\xi)|^2\,d\xi}{(\xi+\kappa+\frac23\vk)^{1-2s}},
\end{align}
uniformly for $\vk\geq 0$.

Consider first the case $\vk=0$ of \eqref{8:star}.  From \eqref{H^s_k defn} and our choice of $\kappa$, we have
\begin{align*}
\| u(t) \|_{H^s} \lesssim  \| u(t) \|_{H^s_\kappa } \lesssim \| u(0) \|_{H^s_\kappa }
		\lesssim  \bigl( 1 + C_0 \| u(0)\|_{L^2}^2 \bigr)^s\| u(0)\|_{L^2} + \| u(0) \|_{H^s},
\end{align*}
which proves \eqref{8:APB}.

When $0\leq s<\frac12$, elementary manipulations show that a bounded subset $U$ of $H^s(\R)$ is $H^s$-equicontinuous if and only if 
$$
\lim_{\vk\to\infty} \ \sup_{u\in U} \ \int_0^\infty \frac{(\xi+\kappa)|\widehat u(\xi)|^2\,d\xi}{(\xi+\kappa+\vk)^{1-2s}} = 0 \qtq{for some} \kappa\geq 0.
$$
In this way, the equicontinuity of $U$ is transferred to $U^*$ by \eqref{8:star}.
\end{proof}

\begin{proof}[Proof of Theorem~\ref{Th: Orbital Stability} for $0\leq s <\frac12$.]
It suffices to prove that for any sequence of solutions $u_n(t)$ satisfying
\begin{align}\label{8:quest}
\inf_{\vec{c}\in\R^N}\, \bigl\|u_{n}(0)-Q_{\Lambda,\vec{c}}\, \bigr\|_{H^s}\to 0, \qtq{we have}
\inf_{\vec{c} \in \mathbb{R}^N} \, \bigl\| u_n(t_n) - Q_{\Lambda, \vec{c}} \,\bigr\|_{H^s} \to 0
\end{align}
for any sequence of times $t_n$.

In the previous section, we showed that Theorem~\ref{Th: Orbital Stability} is valid in $H^{-1/4}(\R)$ and consequently, we do know that for a suitable choice of $\vec c_n$,
\begin{align}\label{8:yay}
\bigl\| u_n(t_n) - Q_{\Lambda, \vec{c}_n} \bigr\|_{H^{-1/4}} \to 0.
\end{align}
To upgrade this convergence to convergence in $H^s(\R)$, it suffices to show $H^s$-equicontinuity of the set $\{u_n(t_n) - Q_{\Lambda, \vec{c}_n}: n\in \N \}$.

Recall from \eqref{Q decay} that $\{Q_{\Lambda, \vec{c}} : \vec c \in \R^N\}$ is $H^1$-bounded; this implies that this set is $H^s$-bounded and equicontinuous.
As the sequence $u_n(0)$ converges in $H^s(\R)$ to this manifold, so $U:=\{u_n(0)\}_{n\geq 1}$ is $H^s$-bounded and equicontinuous.  By Proposition~\ref{P:s pos}, we infer that $\{u_n(t_n)\}_{n\geq 1}$ is also $H^s$-bounded and equicontinuous.

These observations prove that $\{u_n(t_n) - Q_{\Lambda, \vec{c}_n} : n\in \N\}$ is $H^s$-equicontinuous and so \eqref{8:yay} implies convergence in $H^s(\R)$.\end{proof}

Lastly, it remains to treat the energy space $H^{1/2}(\R)$.  This name reflects the fact that it is the natural space on which both the momentum
\begin{equation}\label{Mass}
M(u) := \int   \tfrac{1}{2} |u|^2 \, dx = \langle u_+, u_+\rangle
\end{equation}
and the the Hamiltonian
\begin{equation}\label{Hamil}
E(u) := \int   {\tfrac{1}{2} u H u' } - \tfrac{1}{3} u^3 \, dx = \langle u_+, L_u u_+\rangle
\end{equation}
are well-defined.  (Recall $H$ denotes the Hilbert transform.)  Both the momentum and energy of multisolitons are easily evaluated using Theorem~\ref{T: Sun} and \eqref{Wu relation Hs} (or alternatively, Theorem~\ref{Th: variational problem} and \eqref{beta expand}):
\begin{equation}\label{E(Q)}
M(Q_{\Lambda, \vec{c}}) = \sum_{\lambda\in\Lambda} 2\pi|\lambda|  \qtq{and} E(Q_{\Lambda, \vec{c}} ) = \sum_{\lambda\in\Lambda} -2\pi|\lambda|^2 \, .
\end{equation}

\begin{proof}[Proof of Theorem~\ref{Th: Orbital Stability} for $s=\frac12$.]
Let $u_n$ be a sequence of solutions with
\begin{align}\label{1999}
\inf_{\vec{c}\in\R^N}\, \bigl\|u_{n}(0)-Q_{\Lambda,\vec{c}}\, \bigr\|_{H^{1/2}}\to 0 .
\end{align}

Applying the $s=1/6$ case of Theorem~\ref{Th: Orbital Stability} and Sobolev embedding, we see that for any sequence of times $t_n$ one can find corresponding parameters $\vec c_n$ so that
\begin{align}\label{2000}
\bigl\|u_{n}(t_n)-Q_{\Lambda,\vec{c}_n}\, \bigr\|_{L^3(\R)}\to 0 .
\end{align}
From \eqref{8:APB} and \eqref{Q decay} we know that
\begin{align}\label{2000.5}
\sup_n \| u_n(t_n) \|_{H^{1/2}} + \sup_n  \| Q_{\Lambda,\vec{c}_n} \|_{H^{1/2}} + \sup_n  \| Q_{\Lambda,\vec{c}_n} \|_{W^{1,{3/2}}} < \infty.
\end{align}
By Sobolev embedding, \eqref{2000}, and \eqref{2000.5},
\begin{align}\label{2001}
\Bigl|\int u_n(t_n,x)^3 \, dx & -  \int Q_{\Lambda,\vec c_n}(,x)^3 \, dx\Bigr|\notag\\
&\lesssim \|u_n(t_n)- Q_{\Lambda,\vec c_n}\|_{L^3}\bigl[ \| u_n(t_n) \|_{H^{1/2}}^2 + \| Q_{\Lambda,\vec{c}_n} \|_{H^{1/2}}^2\bigr] \longrightarrow 0
\end{align}
and
\begin{align}\label{2001.5}
\bigl| \bigl\langle u_{n}(t_n)-Q_{\Lambda,\vec{c}_n}, Q_{\Lambda,\vec{c}_n} \bigr\rangle_{H^{1/2}} \bigr|
	\leq \bigl\| u_{n}(t_n)-Q_{\Lambda,\vec{c}_n}\bigr\|_{L^3}  \bigl\| Q_{\Lambda,\vec{c}_n}\bigr\|_{W^{1,3/2}} \longrightarrow 0
\end{align}
as $n\to \infty$.

Returning to \eqref{1999} and using that both $E(u)$ and $M(u)$ are $H^{1/2}$-continuous, conserved under the flow, and constant across the multisoliton manifold, we have
\begin{align}\label{2002}
\Bigl[ E\bigl(u_n(t_n)\bigr)+M\bigl(u_n(t_n)\bigr) \Bigr] \! -\! \Bigl[ E\bigl(Q_{\Lambda,\vec c_n}\bigr)+M\bigl(Q_{\Lambda,\vec c_n}\bigr) \Bigr] \longrightarrow 0 \qtq{as} n\to \infty.
\end{align}
Combining \eqref{2001} and \eqref{2002}, we deduce that
\begin{align*}
\bigl\| u_n(t_n) \bigr\|_{H^{1/2}}^2 -  \bigl\| Q_{\Lambda,\vec c_n} \bigr\|_{H^{1/2}}^2 \longrightarrow 0 \qtq{as} n\to \infty.
\end{align*}
Combining this with \eqref{2001.5}, we find
\begin{align*}
\| u_{n}(t_n)-Q_{\Lambda,\vec{c}_n}\bigr\|_{H^{1/2}}^2&= \| u_n(t_n) \|_{H^{1/2}}^2 - \| Q_{\Lambda,\vec{c}_n} \|_{H^{1/2}}^2\\
&\quad - 2\Re\bigl\langle u_{n}(t_n)-Q_{\Lambda,\vec{c}_n}, Q_{\Lambda,\vec{c}_n} \bigr\rangle_{H^{1/2}}\longrightarrow 0 \qtq{as} n\to \infty,
\end{align*}
which completes the proof of Theorem~\ref{Th: Orbital Stability} when $s=\frac12$.
\end{proof}

\subsection*{Acknowledgements} R.~B. was supported by an AMS-Simons Travel Grant.  R.~K. was supported by NSF grants DMS-2154022, DMS-2452346, and the project ANR-25-CFFS-0004 ``PhysMathEDPInteg'' of the France 2030 program. M.~V. was supported by NSF grant DMS-2348018 and the Simons Fellowship SFI-MPS-SFM-00006244. We thank Gong Chen for discussions at an early stage of this project.  We are also grateful to the referees for their close attention to this paper.

\subsection*{Declarations} No new data was created during this study.  The authors have no financial or proprietary interests in any material discussed in this article.

\bibliographystyle{amsplain}
\bibliography{BOrefs}

\end{document}